\documentclass[review]{elsarticle}

\usepackage{hyperref}
\usepackage{ amsfonts, amssymb}
\usepackage{amsmath}
\usepackage[dvipsnames]{xcolor}
\journal{Journal of \LaTeX\ Templates}
\biboptions{sort&compress}

\newcommand{\dl}{\delta}

\newcommand{\e}{\varepsilon}

\newcommand{\lm}{\lambda}

\newcommand{\s}{\sigma}

\newcommand{\tr}{\operatorname{tr\,}}

\newcommand{\Ree}{\operatorname{Re\,}}

\newcommand{\T}{\operatorname{T\,}}

\newcommand{\diag}{\operatorname{diag\,}}
\newcommand{\wt}{\widetilde}
\begin{document}

\begin{frontmatter}

\title{Asymptotic stability conditions for linear coupled impulsive
systems with time-invariant subsystems}
\tnotetext[mytitlenote]{This research was supported by the German Research Foundation (DFG) via grants No. DA 767/12-1 and SL 343/1-1, and by the German academic exchange service (DAAD), personal ref. no.: 91775148.}

\author[mymainaddress]{Vitalii Slynko\corref{mycorrespondingauthor}}
\cortext[mycorrespondingauthor]{Vitalii Slynko}
\ead{vitalii.slynko@mathematik.uni-wuerzburg.de}

\author[mymainaddress]{Sergey Dashkovskiy}
\ead{sergey.dashkovskiy@mathematik.uni-wuerzburg.de}

\author[mymainaddress]{Ivan Atamas}
\ead{ivan.atamas@mathematik.uni-wuerzburg.de}

\address[mymainaddress]{97074 Emil-Fischer-Str. 40, W\"urzburg}

\begin{abstract}
This article proposes an approach to construct a Lyapunov function for a linear coupled impulsive system consisting of two time-invariant subsystems. In contrast to various variants of small-gain stability conditions for coupled systems, the asymptotic stability property of independent subsystems is not assumed. To analyze the asymptotic stability of a coupled system, the direct Lyapunov method is used in combination with the discretization method. The periodic case and the case when the Floquet theory is not applicable are considered separately. The main results are illustrated with examples.
\end{abstract}

\begin{keyword}
Linear impulsive systems \sep Lyapunov functions\sep Lyapunov stability \sep coupled systems \sep time-variant systems.
\end{keyword}
\end{frontmatter}


\section{Introduction}

Impulsive differential equations can model mechanical systems subjected to shocks. Instantaneous change in holonomic or nonholonomic constraints imposed on the system and changes in parameters of the system in time lead to study impulsive systems with variable coefficients. In this case, it is important to obtain stability conditions that are robust with respect to  variations in the sequence of moments of impulse action.

The theory of impulsive differential equations emerged as an independent direction in the modern differential equations theory and system theory, starting with the classic book \cite{sam-per}.
The main methods for studying stability for systems of differential equations with impulsive action are laid down in \cite{sam-per,LBS}. In modern control theory, impulsive systems are often considered as an important subclass of hybrid systems
\cite{teel,libernzon2}. The stability of linear impulsive systems of differential equations with constant parameters has been a subject of research in many works. In contrast to linear time-invariant systems of ordinary differential equations, where the stability problem is exhaustively solved by the classical Routh-Hurwitz theorem, the stability question is open in the general case. This is due to the fact that the dynamics of impulsive systems is determined not only by the parameters of the system, but also by the sequence of moments of impulse action. It should also be taken into account that an impulsive system can be asymptotically stable even in the case when the continuous and discrete dynamics are both unstable. Therefore, the stability conditions for linear impulsive systems should cover this case as well, see \cite{dv-sl-1,dv-sl-2,briat,briat-1,liu,sl-biv-tunc,sl-biv-1}.  In \cite{sl-biv-tunc,sl-biv-1}, linear impulsive differential equations in Banach spaces were considered under the assumption that the moments of impulse action satisfy the ADT condition, assuming that the continuous and discrete dynamics of the system are both unstable. Using the identities of the commutator calculus, conditions for asymptotic stability are obtained. In \cite{briat,briat-1,BS2015}, the dwell-time conditions which guarantee the asymptotic stability of linear impulsive systems with constant parameters are obtained. In this case, a construction of a Lyapunov function is reduced to an approximate solution of a two-point boundary value problem with boundary conditions in the form of matrix inequalities. A generalization of these results for some classes of linear impulsive systems with variable coefficients is given in \cite{briat-1}. In \cite{dv-sl-2}, the dwell-times estimates that guarantee the asymptotic stability of linear impulsive systems are obtained based on the second Lyapunov method using the derivatives of the second \cite{liu} and higher orders \cite{dv-sl-2} of a Lyapunov function.

The papers \cite{dv-sl-1,liu, dv-sl-3, hadad, dash-miron} are devoted to the study of stability of large-scale impulsive systems. In \cite{dv-sl-1,liu, hadad}, Lyapunov vector functions are used to study the asymptotic stability of the equilibrium of nonlinear impulsive systems. In \cite{dv-sl-3}, the problem of stability of critical equilibria of nonlinear large-scale impulsive systems is considered. In \cite{dash-miron}, the input to state stability theory (ISS) is developed for nonlinear impulsive systems and sufficient conditions for the global asymptotic stability of nonlinear coupled impulsive systems are established (small-gain theorems).
It should be noted that all of these results use the assumption of asymptotic stability of independent subsystems. Hence, it is of interest to find stability conditions for large-scale (coupled) systems which do not assume a priori the presence of the asymptotic stability property of independent subsystems. One of the possible approaches to solve this problem is to use the construction of a matrix-valued Lyapunov function \cite{dj,mart1}. The question of the method of construction is decisive for the practical application of matrix-valued Lyapunov functions in stability problems. In  \cite{mart-sl}, for linear time-variant coupled systems with time-invariant subsystems, a method of construction of a matrix-valued Lyapunov function by some heuristic simplifications of a matrix differential equation written in the block form is proposed. Although, this method allows in some cases to obtain a conclusion about the stability of linear coupled systems with unstable subsystems, the question about the degree of conservatism of the obtained sufficient stability conditions is open. For linear large-scale impulsive systems with variable coefficients, the problem of choosing the elements of a matrix-valued Lyapunov function has not been studied.

Note that for time-discrete coupled systems, some approaches to study Lyapunov stability and ISS in case when some of the independent subsystems do not have the exponential stability property or, respectively, the ISS property are presented in \cite{NWZ2018, NWZ2018(1),GM2017}. For continuous-time and impulsive coupled systems with possibly unstable subsystems, the problem of construction of a Lyapunov function remains open even in the linear case of time-variant systems.

The development of computer calculations has opened up new possibilities to solve the problem of construction of Lyapunov functions for various classes of dynamical systems. In recent years, the use of the discretization method to construct approximate solutions of a Lyapunov matrix differential equations (or their modifications) has led to significant advances in the theory of stability of linear hybrid systems with constant parameters \cite{liron}, systems with delay \cite{xiang} and others. The application of the discretization method of construction of Lyapunov functions makes it possible to obtain with high accuracy estimates of the dwell-times that guarantee the stability of linear hybrid systems or the conditions of robust stability. In \cite{sl-biv-tunc1}, the discretization method is used to synthesize robust control of a nonlinear affine system. \textcolor{blue}{The asymptotic stability of a large-scale system, using the direct Lyapunov method in combination with the discretization method and identities of the commutator calculus were considered in \cite{STA2022}. Here the independent subsystems may not have the asymptotic stability property. In contrast to \cite{STA2022}, we consider a coupled impulsive system by substantially modifying the choice of a candidate of a Lyapunov function using the  time-invariance of (disconnected) independent subsystems.}

In our paper, for the first time, it is proposed to apply the discretization method in order to construct matrix-valued Lyapunov functions for linear impulsive systems with periodic coefficients. In this case, it is assumed that the independent subsystems are time-invariant, and for the dwell-times, two possible cases are considered: they are constant or subject to two-sided estimates. The elements of a matrix-valued Lyapunov function are constructed in the bilinear forms with time-variable matrices. The proposed algorithm of construction of Lyapunov functions admits a simple numerical implementation.

Contributions of this manuscript are as follows.
First of all, a new method for construction of a Lyapunov function for a linear time-variant system consisting of two coupled time-invariant subsystems is proposed.
Conditions for the asymptotic stability of a linear time-variant system are obtained under various assumptions about the dwell-times and dynamic properties of independent subsystems.
 We show that the proposed approach is applicable in the case when one of the subsystems is not stable and the classical methods for studying of coupled systems, which are based on the ideas of  Lyapunov vector functions or small-gain theorems, obtained on the basis of the ISS concept, are not applicable. In the case when subsystems are asymptotically stable, but the small-gain conditions are not satisfy, the proposed approach can still work and leads to less conservative stability conditions. The proposed results are new not only in the context of the theory of impulsive systems, but also for coupled periodic ODEs.

 We introduce examples where the obtained sufficient stability conditions are applicable under different assumptions regarding the continuous and discrete dynamics of the system, in particular, when both dynamics are simultaneously unstable and the application of known results of the theory of stability of impulsive systems is impossible or very difficult.

Next section describes the problem statement. Section 3 is devoted to an informal discussion of the proposed method of construction of a matrix-valued Lyapunov function. In the fourth section, this method is rigorously justified for the case of linear impulsive periodic systems. Sufficient conditions for asymptotic stability are obtained. Section 5 is devoted to the substantiation of the proposed algorithm of construction of a matrix-valued Lyapunov function in the case when the dwell-times are not constant (in this case, the Floquet theory is not applicable). In Section 6 we consider particular cases of rapidly changing interaction between subsystems. We compare our results with known small-gain theorems. Section 7 provides numerical examples that are discussed in Section~8.

{\it Notation.} Let $\Bbb R^n$ be the Euclidian space with standard dot product and $\Bbb R^{n\times m}$ be the linear space of $n\times m$ matrices. For $A\in\Bbb R^{n\times n}$, $\s(A)$ denotes its spectrum, $r_{\s}(A)$ denotes its spectral radius and norm $\|A\|=\lm_{\max}^{1/2}(A^{\T}A)$. If $\sigma(A)\subset\Bbb R$, then $\lm_{\min}(A)$ and $\lm_{\max}(A)$ are its smallest and largest eigenvalues respectively, $\lm_{\max}^+(A)=\max(\lm_{\max}(A),0)$. For any symmetric matrices $P$ and $Q$ notation $P\geq Q$ means that $P-Q$ is a positive semidefinite matrix and $P\succ Q$ means that $P-Q$ is a positive definite matrix. We will also use the Cauchy--Bunyakovsky inequality $|x^{\T}y|\le \|x\|\|y\|$ for $x,y\in\Bbb R^n$.

\section{Statement of the problem} Consider a coupled linear system of impulsive differential equations
\begin{equation}\label{3.1}
\gathered
\dot x_1(t)=A_{11}x_1(t)+A_{12}(t)x_2(t),\quad t\ne\tau_k\\
\dot x_2(t)=A_{21}(t)x_1(t)+A_{22}x_2(t),\quad t\ne\tau_k,\\
x_1(t^+)=B_{11}x_1(t)+B_{12}x_2(t),\quad t=\tau_k,\\
x_2(t^+)=B_{21}x_1(t)+B_{22}x_2(t),\quad t=\tau_k
\endgathered
\end{equation}
where $x_i\in\Bbb R^{n_i}$, $i=1,2$, $A_{ij}\,:\Bbb R\to \Bbb R^{n_i\times n_j}$   are piece-wise continuous functions, $i,j=1,2$.
Suppose that $A_{ii}$ are constant, $A_{ij}(t)$, $i\ne j$ are $\theta$-periodic, i.e., $A_{ij}(t+\theta)=A_{ij}(t)$ for all $t\in\Bbb R$,
$\{\tau_k\}_{k=0}^{\infty}$ is a sequence of moments of impulse action, and
$B_{ij}\in\Bbb R^{n_i\times n_j}$ are constant.
We denote $x=(x_1^{\T},x_2^{\T})^{\T}$ and $n=n_1+n_2$.
For $\{\tau_k\}_{k=0}^{\infty}$, we assume that there are
positive constants $\theta_1$ and $\theta_2$ such that the dwell time $T_k=\tau_k-\tau_{k-1}$, $k\geq 1$ satisfy $\theta_1\le T_k\le\theta_2$.

Coupled systems of the form \eqref{3.1} naturally arise in the process of mathematical modeling of the dynamics of coupled impulsive time-independent systems that exchange information with each other. We consider the case when the coupling functions between subsystems change in time, and independent subsystems are not subject to parametric disturbances and are described by linear systems with constant parameters. Classical approaches to the study of the problem of stability of coupled systems based on a Lyapunov vector function or the concept of ISS a priori assume the property of asymptotic stability of independent subsystems. Therefore, it is of interest to get rid of this a priori assumption, which is caused not by the essence of the problem, but by the restrictions of existing methods for studying stability.

The results presented below can be extended easily to the case of an arbitrary number of independent subsystems. Here we restrict ourselves to the case of two subsystems in order to make the presentation more accessible without overshadowing it with technical details. We study a class of linear non-autonomous systems that allow decomposition into subsystems with time-invariant independent subsystems. The construction of a Lyapunov function for this class can be significantly simplified due to decomposition, in comparison with the direct study of the system \eqref{3.1} without resorting to decomposition.

Note that there exist real constants $\mu_i$, $\delta_i$, $M_i$, $N_i$, $i=1,2$, such that
\begin{equation*}
\gathered
\|e^{sA_{ii}}\|\le M_ie^{s\mu_i},\quad \|e^{-sA_{ii}}\|\le N_ie^{s\delta_i},\quad s\ge 0.
\endgathered
\end{equation*}

{\it Remark 4.1.} There are several ways to obtain the estimate $\|e^{tA}\|\le Me^{\mu t}$, $t\ge 0$, $A\in\Bbb R^{n\times n}$. Here we use a result from \cite{Gil1993}, where one can find the estimate

\begin{equation}\label{matr_exp}
\gathered
\|e^{tA}\|\le e^{\beta_At}\sum\limits_{k=0}^{n-1}\frac{g^k_At^k}{(k!)^{3/2}},
\endgathered
\end{equation}
 where $\beta_A=\max\{\Ree \lambda\,|\lambda\in\sigma(A)\}$, $g_A=\sqrt{\tr(AA^{\T})-|\tr A^2|}$. 
  From \eqref{matr_exp}, we can derive the estimates we need as follows.

If $g_A=0$, $M=1$ and $\mu=\beta_A$ then from \eqref{matr_exp} immediately implies the estimates we need.
Let $g_A\ne 0$ and for a given $\epsilon>0$, we denote
\begin{equation*}
\gathered
M_{\epsilon,A}:=\sup\limits_{t\ge 0}e^{-\epsilon t}\sum\limits_{k=0}^{n-1}\frac{g^k_At^k}{(k!)^{3/2}}<\infty,
\endgathered
\end{equation*}
then 
\begin{equation*}
\gathered
\|e^{tA}\|\le M_{\epsilon,A}e^{(\beta_A+\epsilon)t}, \quad t\ge 0.
\endgathered
\end{equation*}

Since the functions $A_{ij}(t)$ are assumed to be bounded and periodic, then for some positive constants $\gamma_{12}^{(m)}$, $\gamma_{21}^{(m)}$, $m=0,1,\dots,N-1$ it holds that
\begin{equation*}
\gathered
\sup\limits_{s\in(mh,(m+1)h]}\|A_{12}(s)\|\le\gamma_{12}^{(m)},\quad
\sup\limits_{s\in(mh,(m+1)h]}\|A_{21}(s)\|\le\gamma_{21}^{(m)}.
\endgathered
\end{equation*}
For any $m\in\Bbb Z$, let $\gamma_{ij}^{(m)}:=\gamma_{ij}^{(\varrho)}$, where $\varrho$ is the remainder of $m$ devided by $N$.
By this we extend the definition of constants $\gamma_{ij}^{(m)}$ for any $m\in\Bbb Z$.
Here we study the asymptotic stability of \eqref{3.1} in the sense of:

{\bf Definition 2.1.} System of differential equations \eqref{3.1} is called

1) stable if for any $\e>0$, $t_0\in\Bbb R$ there exists
$\dl=\dl(\e,t_0)>0$ such that the inequality $\|x_0\|<\dl$ implies that $\|x(t,t_0,x_0)\|<\e$ for all $t\ge t_0$;

2) uniformly stable if for any $\e>0$ there exists
$\dl=\dl(\e)>0$ such that for all $t_0\in\Bbb R$ the inequality $\|x_0\|<\dl$ implies  $\|x(t,t_0,x_0)\|<\e$ for $t\ge t_0$;

3) asymptotically stable if it is stable and
$\lim\limits_{t\to+\infty}\|x(t,t_0,x_0)\|=0$.

Here, $x(t,t_0,x_0)$ is the solution of the Cauchy problem \eqref{3.1} with the initial condition $x(t_0,t_0,x_0)=x_0$, $x_0=(x_{10}^{\T},x_{20}^{\T})^{\T}\in\Bbb R^n$.

The aim of this work is to construct a Lyapunov function for \eqref{3.1} and to prove sufficient conditions for its stability.
The case of the periodic system \eqref{3.1}, when $\theta_1=\theta_2=\theta$ and the general case, when
$\theta_1<\theta_2$ are considered separately. If the linear impulsive system \eqref{3.1} is not periodic, the Floquet theory is not applicable.

\section{General idea of construction of a matrix-valued Lyapunov function} 

\textcolor{blue}{Here we consider the case $\theta_1=\theta_2$ and present the general idea, how we will derive a Lyapunov function for \eqref{3.1}. Without any restriction we assume that $t_0=0$. Let $\mathfrak V(t,x)=(v_{ij}(t,\cdot,\cdot))_{i,j=1,2}$ be a matrix-valued Lyapunov function (MFL) \cite{dj,mart2}. We choose $v_{ij}(t,x_i,x_j)=x_i^{\T}P_{ij}(t)x_j$, where $P_{ij}\,:\Bbb R\to\Bbb R^{n_i\times n_j}$ are continuous on the left $\theta$-periodic maps, $P_{ij}(t)=P_{ji}^{\T}(t)$. It is enough to define $P_{ij}(t)$, $i,j=1,2$ on the period $(0,\theta]$.
For the practical application of this rather general design, we provide a method of construction of $P_{ij}(t)$.}

\textcolor{blue}{From $\mathfrak V$ we proceed to the following scalar Lyapunov function
\begin{equation*}
\gathered
v(t,x_1,x_2)=(1,1)\mathfrak V(t,x_1,x_2) (1,1)^{\T}=v_{11}(t,x_1)+2v_{12}(t,x_1,x_2)+v_{22}(t,x_2)\\
=(x_1^{\T},x_2^{\T})P(t)(x_1^{\T},x_2^{\T})^{\T},
\endgathered
\end{equation*}
where $P(t)=(P_{ij}(t))_{i,j=1,2}$ is a block matrix $P_{ij}(t)=P_{ji}^{\T}(t)$ to be designed.}

\textcolor{blue}{Let matrix $P(t)$ satisfy the condition
\begin{equation}\label{4.1}
\gathered
\dot P(t)+A^{\T}(t)P(t)+P(t)A(t)=0,\quad t\in(k\theta,(k+1)\theta),\quad k\in\Bbb Z_+,
\endgathered
\end{equation}
with the block matrix $A(t)=(A_{ij}(t))_{i,j=1,2}$. It is easy to show that for $t=k\theta$ the following estimate holds
\begin{equation*}
\gathered
v(k\theta+0,x_1(k\theta+0),x_2(k\theta+0))=x^{\T}(k\theta+0)P(k\theta+0)x(k\theta+0)\\=
(Bx(k\theta))^{\T}P_0Bx(k\theta)=x^{\T}(k\theta)B^{\T}P_0Bx(k\theta)
\le\lm v(k\theta,x_1(k\theta),x_2(k\theta)),
\endgathered
\end{equation*}
where $\lm=\lm_{\max}(B^{\T}P_0B(P(\theta))^{-1})$, $B=(B_{ij})_{i,j=1,2}$ and $P_0=P(0+0)$ is a symmetric positive definite matrix. From \eqref{4.1} it follows that $v(t,x_1,x_2)$ satisfies the following system of impulsive differential inequalities
\begin{equation*}
\gathered
\dot v(t,x_1(t),x_2(t))=0,\quad t\ne k\theta,\\
v(t+0,x_1(t+0),x_2(t+0))\le\lm v(t,x_1(t),x_2(t)),\quad t=k\theta.
\endgathered
\end{equation*}
By means of the comparison principle \cite{LBS} the stability investigation of \eqref{3.1} reduces to the stability question of 
\begin{equation}\label{LCS}
\gathered
\dot u(t)=0,\quad t\ne k\theta,\\
u(t+0)=\lm u(t),\quad t=k\theta.
\endgathered
\end{equation}
The asymptotic stability of \eqref{LCS} implies the same property for \eqref{3.1}. The condition $\lm<1$ (which restricts $P_0\succ 0$) is necessary and sufficient for the asymptotic stability of \eqref{LCS} and is equivalent to the necessary and sufficient condition for the asymptotic stability of \eqref{3.1},
as can be seen from the Floquet-Lyapunov theorem. Resolving \eqref{4.1} on the interval $(k\theta,(k+1)\theta)$ is equivalent to the calculation of the monodromy matrix $\Phi$ for \eqref{3.1}, i.e., the problem is comparable in complexity to the integration of \eqref{3.1}. The main idea of construction of the matrix-valued Lyapunov function is to construct an approximate solution to \eqref{4.1} written in the block form:
\begin{equation}\label{MEq}
\gathered
\dot P_{11}(t)+A_{11}^{\T}P_{11}(t)+P_{11}(t)A_{11}=-(P_{12}(t)A_{21}(t)+A_{21}^{\T}(t)P_{21}(t)),\\
\dot P_{22}(t)+A_{22}^{\T}P_{22}(t)+P_{22}(t)A_{22}=-(A_{12}^{\T}(t)P_{12}(t)+P_{21}(t)A_{12}(t)),\\
\dot P_{12}(t)+A_{11}^{\T}P_{12}(t)+P_{12}(t)A_{22}=-(P_{11}(t)A_{12}(t)+A_{21}^{\T}(t)P_{22}(t)).
\endgathered
\end{equation}}

\textcolor{blue}{Since $P(t)$ is $\theta$-periodic, we need to construct a solution to \eqref{MEq} on the interval $(0,\theta]$. Let $P_{ij}(0+0)=P_{ij}^{(0)}$, $i,j=1,2$, $P_{ij}^{(0)}=(P_{ji}^{(0)})^{\T}$, where $(P_{ij}^{(0)})_{i,j=1,2}$ is a positive definite block matrix.
We use the discretization method to resolve \eqref{MEq}.
For $N\in\mathbb N$ let $h=\frac{\theta}{N}$ be the discretization step.
We derive an approximate solution to \eqref{MEq} step by step on the intervals $(mh,(m+1)h]$, $m=0,\dots,N-1$.
By the Cauchy formula applied to the interval $(mh,(m+1)h]$ we obtain integral representations for the solutions of \eqref{MEq}
\begin{equation*}
\gathered
P_{11}(t)=e^{-A_{11}^{\T}(t-mh)}P_{11}^{(m)}e^{-A_{11}(t-mh)}\\
-\int\limits_{mh}^te^{-A_{11}^{\T}(t-s)}(P_{12}(s)A_{21}(s)+A_{21}^{\T}(s)P_{21}(s))e^{-A_{11}(t-s)}\,ds,\\
\endgathered
\end{equation*}
\begin{equation*}
\gathered
P_{22}(t)=e^{-A_{22}^{\T}(t-mh)}P_{22}^{(m)}e^{-A_{22}(t-mh)}\\
-\int\limits_{mh}^te^{-A_{22}^{\T}(t-s)}(A_{12}^{\T}(s)P_{12}(s)+P_{21}(s)A_{12}(s))e^{-A_{22}(t-s)}\,ds,\\
\endgathered
\end{equation*}
\begin{equation*}
\gathered
P_{12}(t)=e^{-A_{11}^{\T}(t-mh)}P_{12}^{(m)}e^{-A_{22}(t-mh)}\\
-\int\limits_{mh}^te^{-A_{11}^{\T}(t-s)}(P_{11}(s)A_{12}(s)+A_{21}^{\T}(s)P_{22}(s))e^{-A_{22}^{\T}(t-s)}\,ds.
\endgathered
\end{equation*}
Here, $P_{ij}^{(m)}=P_{ij}(mh-0)$, $m=1,\dots,N-1$, $i,j=1,2$. 
For $h$ small enough we use the approximation $P_{ij}(s)\approx P_{ij}^{(m)}$ in the integrals which leads to the approximate solutions to \eqref{MEq} for $t\in(mh,(m+1)h]$ as follows}

\textcolor{blue}{\begin{equation*}
\gathered
P_{11}(t)\approx e^{-A_{11}^{\T}(t-mh)}(P_{11}^{(m)}
-\int\limits_{mh}^t(P_{12}^{(m)}A_{21}(s)+A_{21}^{\T}(s)P_{21}^{(m)})\,ds)e^{-A_{11}(t-mh)},
\endgathered
\end{equation*}
\begin{equation*}
\gathered
P_{22}(t)\approx e^{-A_{22}^{\T}(t-mh)}(P_{22}^{(m)}
-\int\limits_{mh}^t(A_{12}^{\T}(s)P_{12}^{(m)}+P_{21}^{(m)}A_{12}(s))\,ds)e^{-A_{22}(t-mh)},
\endgathered
\end{equation*}
\begin{equation*}
\gathered
P_{12}(t)\approx e^{-A_{11}^{\T}(t-mh)}(P_{12}^{(m)}
-\int\limits_{mh}^t(P_{11}^{(m)}A_{12}(s)+A_{21}^{\T}(s)P_{22}^{(m)})\,ds)e^{-A_{22}(t-mh)}.
\endgathered
\end{equation*}}

\textcolor{blue}{Note that for $h\to 0+$ these approximations converge to the true solutions to \eqref{MEq}. Since the asymptotic stability conditions derived with help of $\mathfrak V(t,x_1,x_2)$ are necessary and sufficient, the matrix-valued Lyapunov function, whose elements are given by these approximations leads to sufficient asymptotic stability conditions for \eqref{3.1} arbitrarily close to the necessary ones.}

\section{The case when the system is periodic}

\subsection{Construction of the matrix-valued Lyapunov function}
We proceed to a rigorous justification of the proposed method of construction in the case when the dwell-times are constant and equal to the period $\theta$ of $A_{ij}(t)$.
We introduce discretization parameters: the number of nodes $N\in\Bbb N$ and the discretization step length $h=\frac{\theta}{N}$.
Let $P_0=(P_{ij}^{(0)})_{i,j=1,2}$ be a positive definite symmetric block matrix, $P_{ij}^{(0)}\in\Bbb R^{n_i\times n_j}$, $P_{ij}^{(0)}=(P_{ji}^{(0)})^{\T}$. We define recursively the matrices
$
P_{ij}^{(m)},\quad P_{ji}^{(m)}=(P_{ij}^{(m)})^{\T},\quad i,j=1,2
$
as follows
\begin{equation}\label{5.1}
\gathered
P_{11}^{(m+1)}= e^{-A_{11}^{\T}h}(P_{11}^{(m)}
-\int\limits_{mh}^{(m+1)h}(P_{12}^{(m)}A_{21}(s)+A_{21}^{\T}(s)P_{21}^{(m)})\,ds)e^{-A_{11}h},
\endgathered
\end{equation}
\begin{equation}\label{5.2}
\gathered
P_{22}^{(m+1)}=e^{-A_{22}^{\T}h}(P_{22}^{(m)}
-\int\limits_{mh}^{(m+1)h}(A_{12}^{\T}(s)P_{12}^{(m)}+P_{21}^{(m)}A_{12}(s))\,ds)e^{-A_{22}h}
\endgathered
\end{equation}
\begin{equation}\label{5.3}
\gathered
P_{12}^{(m+1)}= e^{-A_{11}^{\T}h}(P_{12}^{(m)}
-\int\limits_{mh}^{(m+1)h}(P_{11}^{(m)}A_{12}(s)+A_{21}^{\T}(s)P_{22}^{(m)})\,ds)e^{-A_{22}h}.
\endgathered
\end{equation}
Next, define the matrices $P_{ij}(t)$, $i,j=1,2$, $P_{ij}(t)=P_{ji}^{\T}(t)$ on the intervals $(mh,(m+1)h]$ by setting
\begin{equation}\label{5.4}
\gathered
P_{11}(t)= e^{-A_{11}^{\T}(t-mh)}(P_{11}^{(m)}
-\int\limits_{mh}^t(P_{12}^{(m)}A_{21}(s)+A_{21}^{\T}(s)P_{21}^{(m)})\,ds)e^{-A_{11}(t-mh)},
\endgathered
\end{equation}
\begin{equation}\label{5.5}
\gathered
P_{22}(t)=e^{-A_{22}^{\T}(t-mh)}(P_{22}^{(m)}
-\int\limits_{mh}^t(A_{12}^{\T}(s)P_{12}^{(m)}+P_{21}^{(m)}A_{12}(s))\,ds)e^{-A_{22}(t-mh)}
\endgathered
\end{equation}
\begin{equation}\label{5.6}
\gathered
P_{12}(t)= e^{-A_{11}^{\T}(t-mh)}(P_{12}^{(m)}
-\int\limits_{mh}^t(P_{11}^{(m)}A_{12}(s)+A_{21}^{\T}(s)P_{22}^{(m)})\,ds)e^{-A_{22}(t-mh)}.
\endgathered
\end{equation}
We define $\mathfrak V(t,x_1,x_2)=(v_{ij}(t,.,.))_{i,j=1,2}$, where  $\,v_{ij}(t,x_i,x_j)=x_i^{\T}P_{ij}(t)x_j$, $P_{ij}(t)=P_{ji}^{\T}(t)$, $i,j=1,2$.
Using the function $\mathfrak V$, we construct the scalar Lyapunov function \cite{dj}
\begin{equation}\label{5.7}
\gathered
v(t,x_1,x_2)=v_{11}(t,x_1)+2v_{12}(t,x_1,x_2)+v_{22}(t,x_2).
\endgathered
\end{equation}
We establish some auxiliary estimates for the derivatives of the components of $\mathfrak V$ along the solutions of \eqref{3.1} and estimates for the Lyapunov function $v(t,x_1,x_2)$, necessary for construction of a comparison linear impulsive differential equation  and obtaining sufficient conditions for the asymptotic stability of  \eqref{3.1}. We will also use the following proposition.

{\bf Proposition 4.1} Let $M,\mu>0$ be such that
$\|e^{tA}\|\le Me^{t\mu}$, $t\ge 0$. Then
\begin{equation}\label{5.8***}
\gathered
\|e^{tA}-I\|\le \frac{\|A\|M}{\mu}(e^{t\mu}-1), \quad t\geq 0.
\endgathered
\end{equation}

{\it Proof.}
Let $X(t)=e^{tA}-I$, then $X(0)=0$
and $\dot X(t)=A(X(t)+I)$, applying the Cauchy formula, we obtain
$X(t)=\int\limits_0^te^{A(t-s)}A\,ds$, which implies \eqref{5.8***}.


The following assertion is necessary to verify the positive-definiteness conditions for the proposed Lyapunov function $v(t,x_1,x_2)$. Since this function depends explicitly on time $t$, it is impossible to check this condition pointwise for all $t \in [0,\theta]$. The following Lemma 4.1 reduces this checking to a finite number of conditions.

{\bf Lemma 4.1.} Let $z_{1m}=e^{-A_{11}(t-mh)}x_1$, $z_{2m}=e^{-A_{22}(t-mh)}x_2$, $t\in(mh,(m+1)h]$. Then,
\begin{equation}\label{5.8}
\gathered
\lm_{\min}(\Pi_m)\|z_m\|^2\le v(t,x_1,x_2)\le\lm_{\max}(\Xi_m)\|z_m\|^2,\\
\text{for all} \quad t\in(mh,(m+1)h],\quad m=0,\dots,N-1,
\endgathered
\end{equation}
where $z_m=(z_{1m}^{\T},z_{2m}^{\T})^{\T}$, $\|z_m\|^2=\|z_{1m}\|^2+\|z_{2m}\|^2$, $\Pi_m=(\pi_{ij}^{(m)})_{i,j=1,2}$, $\Xi_m=(\xi_{ij}^{(m)})_{i,j=1,2}$ are block matrices with the elements
\begin{equation}\label{l4.1}
\gathered
\pi_{11}^{(m)}=P_{11}^{(m)}-h(2\gamma_{21}^{(m)}\|P_{12}^{(m)}\|+(\gamma_{12}^{(m)}\|P_{11}^{(m)}\|+\gamma_{21}^{(m)}\|P_{22}^{(m)}\|))I_{n_1},\\
\pi_{22}^{(m)}=P_{22}^{(m)}-h(2\gamma_{12}^{(m)}\|P_{12}^{(m)}\|+(\gamma_{12}^{(m)}\|P_{11}^{(m)}\|+\gamma_{21}^{(m)}\|P_{22}^{(m)}\|))I_{n_2},\\
\pi_{12}^{(m)}=P_{12}^{(m)},\quad \pi_{21}^{(m)}=P_{21}^{(m)}
\endgathered
\end{equation}
\begin{equation*}
\gathered
\xi_{11}^{(m)}=P_{11}^{(m)}+h(2\gamma_{21}^{(m)}\|P_{12}^{(m)}\|+(\gamma_{12}^{(m)}\|P_{11}^{(m)}\|+\gamma_{21}^{(m)}\|P_{22}^{(m)}\|))I_{n_1},\\
\xi_{22}^{(m)}=P_{22}^{(m)}+h(2\gamma_{12}^{(m)}\|P_{12}^{(m)}\|+(\gamma_{12}^{(m)}\|P_{11}^{(m)}\|+\gamma_{21}^{(m)}\|P_{22}^{(m)}\|))I_{n_2},\\
\xi_{12}^{(m)}=P_{12}^{(m)},\quad \xi_{21}^{(m)}=P_{21}^{(m)}.
\endgathered
\end{equation*}

The proof of this statement is given in the Appendix.

Let

\begin{equation*}
\gathered
\eta_{11}^{(m)}:=\sqrt{\|P_{11}^{(m)}\|^2+\|P_{12}^{(m)}\|^2}+\|P_{12}^{(m)}\|,\\
\eta_{22}^{(m)}:=\sqrt{\|P_{22}^{(m)}\|^2+\|P_{12}^{(m)}\|^2}+\|P_{12}^{(m)}\|
\endgathered
\end{equation*}
\begin{equation*}
\gathered
\eta_{12}^{(m)}:=\frac{1}{2}(\|P_{11}^{(m)}\|\gamma_{12}^{(m)}+\|P_{22}^{(m)}\|\gamma_{21}^{(m)}\\
+\sqrt{(\|P_{11}^{(m)}\|\gamma_{12}^{(m)}+\|P_{22}^{(m)}\|\gamma_{21}^{(m)})^2
+16(\gamma_{21}^{(m)})^2\|P_{12}^{(m)}\|^2}).
\endgathered
\end{equation*}
\begin{equation*}
\gathered
\eta_{21}^{(m)}:=\frac{1}{2}(\|P_{11}^{(m)}\|\gamma_{12}^{(m)}+\|P_{22}^{(m)}\|\gamma_{21}^{(m)}\\
+\sqrt{(\|P_{11}^{(m)}\|\gamma_{12}^{(m)}+\|P_{22}^{(m)}\|\gamma_{21}^{(m)})^2
+16(\gamma_{12}^{(m)})^2\|P_{12}^{(m)}\|^2}).
\endgathered
\end{equation*}

We denote
\begin{equation*}
\gathered
\alpha_{11}^{(m)}:=
\gamma_{12}^{(m)}\|A_{22}\|N_1M_2\eta_{11}^{(m)},\quad
\alpha_{12}^{(m)}:=\gamma_{12}^{(m)}\|A_{11}\|N_1\eta_{11}^{(m)},
\endgathered
\end{equation*}
\begin{equation*}
\gathered
\alpha_{21}^{(m)}:=
\gamma_{21}^{(m)}\|A_{11}\|N_2M_1\eta_{22}^{(m)},\quad
\alpha_{22}^{(m)}:=\gamma_{21}^{(m)}\|A_{22}\|N_2\eta_{22}^{(m)},
\endgathered
\end{equation*}

\begin{equation}\label{Theta}
\gathered
\Theta_m(h):=\frac{\alpha_{11}^{(m)}}{\mu_2}\Big(\frac{e^{(\mu_2+\delta_1)h}-1}{\mu_2+\delta_1}-
\frac{e^{\delta_1h}-1}{\delta_1}\Big)
+\frac{\alpha_{12}^{(m)}}{\delta_1}\Big(\frac{e^{h\delta_1}-1}{\delta_1}-h\Big)\\
+\frac{\alpha_{21}^{(m)}}{\mu_1}\Big(\frac{e^{(\mu_1+\delta_2)h}-1}{\mu_1+\delta_2}-
\frac{e^{\delta_2h}-1}{\delta_2}\Big)
+\frac{\alpha_{22}^{(m)}}{\delta_2}\Big(\frac{e^{h\delta_2}-1}{\delta_2}-h\Big)+
\\
+2\Big(\gamma_{12}^{(m)}\eta_{12}^{(m)}N_1M_2\Big(\frac{he^{(\mu_2+\delta_1)h}}{\mu_2+\delta_1}-\frac{e^{(\mu_2+\delta_1)h}-1}{(\mu_2+\delta_1)^2}\Big)\\
+\gamma_{21}^{(m)}\eta_{21}^{(m)}N_2M_1\Big(\frac{he^{(\mu_1+\delta_2)h}}{\mu_1+\delta_2}-\frac{e^{(\mu_1+\delta_2)h}-1}{(\mu_1+\delta_2)^2}\Big)\Big).
\endgathered
\end{equation}

The following lemma establishes estimates for the change of the Lyapunov function on each of the discretization intervals $(mh,(m+1)h]$.

{\bf Lemma 4.2}
If $\Pi_m$ are positive definite, for all $m=0,\dots,N-1$, then

\begin{equation}\label{5.9}
\gathered
v((m+1)h,x_1((m+1)h),x_2((m+1)h))\\
\le e^{\frac{\Theta_m(h)}{\lm_{\min}(\Pi_m)}}v(mh+0,x_1(mh+0),x_2(mh+0)).
\endgathered
\end{equation}
{\it Proof.} Let $t\in (mh, (m+1)h)$. It can be established by direct calculations that
\begin{equation}\label{5.10}
\gathered
\dot v(t,x_1(t),x_2(t))
=x_1^{\T}(t)(P_{12}(t)A_{21}(t)+A_{21}^{\T}(t)P_{21}(t)\\
-e^{-A_{11}^{\T}(t-mh)}(P_{12}^{(m)}A_{21}(t)+A_{21}^{\T}(t)P_{21}^{(m)})e^{-A_{11}(t-mh)})x_1(t)\\+
2x_1^{\T}(t)(P_{11}(t)A_{12}(t)+A_{21}^{\T}(t)P_{22}(t)\\
-e^{-A_{11}^{\T}(t-mh)}(P_{11}^{(m)}A_{12}(t)+A_{21}^{\T}(t)P_{22}^{(m)})e^{-A_{22}(t-mh)})x_2(t)\\
+x_2^{\T}(t)((A_{12}^{\T}(t)P_{12}(t)+P_{21}(t)A_{12}(t))\\
-e^{-A_{22}^{\T}(t-mh)}(A_{12}^{\T}(t)P_{12}^{(m)}+P_{21}^{(m)}A_{12}(t))e^{-A_{22}(t-mh)})x_2(t)\\
=z_{1m}^{\T}(t)(e^{A_{11}^{\T}(t-mh)}(P_{12}(t)A_{21}(t)+A_{21}^{\T}(t)P_{21}(t))e^{A_{11}(t-mh)}\\
-(P_{12}^{(m)}A_{21}(t)+A_{21}^{\T}(t)P_{21}^{(m)}))z_{1m}(t)\\
+2z_{1m}^{\T}(t)(e^{A_{11}^{\T}(t-mh)}(P_{11}(t)A_{12}(t)+A_{21}^{\T}(t)P_{22}(t))e^{A_{22}(t-mh)}\\
-(P_{11}^{(m)}A_{12}(t)+A_{21}^{\T}(t)P_{22}^{(m)}))z_{2m}(t)\\
+z_{2m}^{\T}(t)(e^{A_{22}^{\T}(t-mh)}(A_{12}^{\T}(t)P_{12}(t)+P_{21}(t)A_{12}(t))e^{A_{22}(t-mh)}\\
-(A_{12}^{\T}(t)P_{12}^{(m)}+P_{21}^{(m)}A_{12}(t)))z_{2m}(t)
\endgathered
\end{equation}
Let us consider separately
\begin{equation*}
\gathered
e^{A_{11}^{\T}(t-mh)}(P_{12}(t)A_{21}(t)+A_{21}^{\T}(t)P_{21}(t))e^{A_{11}(t-mh)}
-(P_{12}^{(m)}A_{21}(t)+A_{21}^{\T}(t)P_{21}^{(m)})\\
=e^{A_{11}^{\T}(t-mh)}P_{12}(t)A_{21}(t)e^{A_{11}(t-mh)}-P_{12}^{(m)}A_{21}(t)\\
+e^{A_{11}^{\T}(t-mh)}A_{21}^{\T}(t)P_{21}(t)e^{A_{11}(t-mh)}-A_{21}^{\T}(t)P_{21}^{(m)}
\endgathered
\end{equation*}
Taking into account the explicit expressions for $P_{12}(t)$ from
\eqref{5.6} we get
\begin{equation*}
\gathered
e^{A_{11}^{\T}(t-mh)}P_{12}(t)A_{21}(t)e^{A_{11}(t-mh)}-P_{12}^{(m)}A_{21}(t)\\=
(P_{12}^{(m)}
-\int\limits_{mh}^t(P_{11}^{(m)}A_{12}(s)+A_{21}^{\T}(s)P_{22}^{(m)})\,ds)e^{-A_{22}(t-mh)}A_{21}(t)e^{A_{11}(t-mh)}\\
-P_{12}^{(m)}A_{21}(t))=P_{12}^{(m)}(e^{-A_{22}(t-mh)}A_{21}(t)e^{A_{11}(t-mh)}-A_{21}(t))\\
-\int\limits_{mh}^t(P_{11}^{(m)}A_{12}(s)+A_{21}^{\T}(s)P_{22}^{(m)})\,ds)e^{-A_{22}(t-mh)}A_{21}(t)e^{A_{11}(t-mh)}
\endgathered
\end{equation*}

Consequently, using \eqref{5.8***} we obtain 
\begin{equation*}
\gathered
\|e^{A_{11}^{\T}(t-mh)}P_{12}(t)A_{21}(t)e^{A_{11}(t-mh)}-P_{12}^{(m)}A_{21}(t)\|\\
\le\gamma_{21}^{(m)}\Big(\|P_{12}^{(m)}\|(\frac{\|A_{11}\|M_1N_2}{\mu_1}e^{\delta_2(t-mh)}(e^{\mu_1(t-mh)}-1)
+\frac{N_2\|A_{22}\|}{\delta_2}(e^{\delta_2(t-mh)}-1))\\
+(t-mh)(\|P_{11}^{(m)}\|\gamma_{12}^{(m)}+\gamma_{21}^{(m)}\|P_{22}^{(m)}\|)M_1N_2e^{(t-mh)(\mu_1+\delta_2)}\Big)
\endgathered
\end{equation*}
and 
\begin{equation*}
\gathered
\|e^{A_{11}^{\T}(t-mh)}A_{21}^{\T}(t)P_{21}(t)e^{-A_{11}(t-mh)}-A_{21}^{\T}(t)P_{21}^{(m)}\|\\
\le\gamma_{21}^{(m)}\Big(\|P_{12}^{(m)}\|(\frac{\|A_{11}\|M_1N_2}{\mu_1}e^{\delta_2(t-mh)}(e^{\mu_1(t-mh)}-1)
+\frac{N_2\|A_{22}\|}{\delta_2}(e^{\delta_2(t-mh)}-1))\\
+(t-mh)(\|P_{11}^{(m)}\|\gamma_{12}^{(m)}+\gamma_{21}^{(m)}\|P_{22}^{(m)}\|)M_1N_2e^{(t-mh)(\mu_1+\delta_2)}\Big)
\endgathered
\end{equation*}

Hence, applying the triangle inequality we obtain 
\begin{equation*}
\gathered
\|e^{A_{11}^{\T}(t-mh)}(P_{12}(t)A_{21}(t)+A_{21}^{\T}(t)P_{21}(t))e^{A_{11}(t-mh)}
-(P_{12}^{(m)}A_{21}(t)+A_{21}^{\T}(t)P_{21}^{(m)})\|\\
\le 2\gamma_{21}^{(m)}\Big(\|P_{12}^{(m)}\|(\frac{\|A_{11}\|M_1N_2}{\mu_1}e^{\delta_2(t-mh)}(e^{\mu_1(t-mh)}-1)
+\frac{N_2\|A_{22}\|}{\delta_2}(e^{\delta_2(t-mh)}-1))\\
+(t-mh)(\|P_{11}^{(m)}\|\gamma_{12}^{(m)}+\gamma_{21}^{(m)}\|P_{22}^{(m)}\|)M_1N_2e^{(t-mh)(\mu_1+\delta_2)}\Big):=\psi_{11}^{(m)}(t).
\endgathered
\end{equation*}
Let us transform the following expression
\begin{equation*}
\gathered
e^{A_{11}^{\T}(t-mh)}(P_{11}(t)A_{12}(t)+A_{21}^{\T}(t)P_{22}(t))e^{A_{22}(t-mh)}
-(P_{11}^{(m)}A_{12}(t)+A_{21}^{\T}(t)P_{22}^{(m)}) \\
=e^{A_{11}^{\T}(t-mh)}P_{11}(t)A_{12}(t)e^{A_{22}(t-mh)}-P_{11}^{(m)}A_{12}(t)\\
+e^{A_{11}^{\T}(t-mh)}A_{21}^{\T}(t)P_{22}(t)e^{A_{22}(t-mh)}-A_{21}^{\T}(t)P_{22}^{(m)}
\endgathered
\end{equation*}
Taking into account \eqref{5.4}, consider separately
\begin{equation*}
\gathered
e^{A_{11}^{\T}(t-mh)}P_{11}(t)A_{12}(t)e^{A_{22}(t-mh)}-P_{11}^{(m)}A_{12}(t)\\
=P_{11}^{(m)}(e^{-A_{11}(t-mh)}A_{12}(t)e^{A_{22}(t-mh)}-A_{12}(t))\\
-\int\limits_{mh}^t(P_{12}^{(m)}A_{21}(s)+A_{21}^{\T}(s)P_{21}^{(m)})\,dse^{-A_{11}(t-mh)}A_{12}(t)e^{A_{22}(t-mh)}
\endgathered
\end{equation*}
Consequently, using \eqref{5.8***} we get 
\begin{equation*}
\gathered
\|e^{A_{11}^{\T}(t-mh)}P_{11}(t)A_{12}(t)e^{A_{22}(t-mh)}-P_{11}^{(m)}A_{12}(t)\|\\\le
\gamma_{12}^{(m)}\Big(\|P_{11}^{(m)}\|(\frac{\|A_{22}\|M_2N_1}{\mu_2}e^{\delta_1(t-mh)}(e^{\mu_2(t-mh)}-1)
+\frac{N_1\|A_{11}\|}{\delta_1}(e^{\delta_1(t-mh)}-1))\\
+2(t-mh)\|P_{12}^{(m)}\|\gamma_{21}^{(m)}N_1M_2e^{(t-mh)(\delta_1+\mu_2)}\Big).
\endgathered
\end{equation*}
Taking into account \eqref{5.5} we obtain
\begin{equation*}
\gathered
e^{A_{11}^{\T}(t-mh)}A_{21}^{\T}(t)P_{22}(t)e^{A_{22}(t-mh)}-A_{21}^{\T}(t)P_{22}^{(m)}\\
=(e^{A_{11}^{\T}(t-mh)}A_{21}^{\T}(t)e^{-A_{22}^{\T}(t-mh)}-A_{21}^{\T}(t))P_{22}^{(m)}\\
-e^{A_{11}^{\T}(t-mh)}A_{21}^{\T}(t)e^{-A_{22}^{\T}(t-mh)}\int\limits_{mh}^t(A_{12}^{\T}(s)P_{12}^{(m)}+P_{21}^{(m)}A_{12}(s))\,ds
\endgathered
\end{equation*}
Thus,  taking into account \eqref{5.8***} we have 
\begin{equation*}
\gathered
\|e^{A_{11}^{\T}(t-mh)}A_{21}^{\T}(t)P_{22}(t)e^{A_{22}(t-mh)}-A_{21}^{\T}(t)P_{22}^{(m)}\|\\
\le \gamma_{21}^{(m)}\Big(\|P_{22}^{(m)}\|(\frac{\|A_{11}\|M_1N_2}{\mu_1}e^{\delta_2(t-mh)}(e^{\mu_1(t-mh)}-1)
+\frac{N_2\|A_{22}\|}{\delta_2}(e^{\delta_2(t-mh)}-1))\\
+2(t-mh)\|P_{12}^{(m)}\|\gamma_{12}^{(m)}N_2M_1e^{(t-mh)(\delta_2+\mu_1)}\Big)
\endgathered
\end{equation*}
Thereby, finally we find the following estimate
\begin{equation*}
\gathered
\|e^{A_{11}^{\T}(t-mh)}(P_{11}(t)A_{12}(t)+A_{21}^{\T}(t)P_{22}(t))e^{A_{22}(t-mh)}
-(P_{11}^{(m)}A_{12}(t)+A_{21}^{\T}(t)P_{22}^{(m)})\| \\
\le \gamma_{12}^{(m)}\Big(\|P_{11}^{(m)}\|(\frac{\|A_{22}\|M_2N_1}{\mu_2}e^{\delta_1(t-mh)}(e^{\mu_2(t-mh)}-1)
+\frac{N_1\|A_{11}\|}{\delta_1}(e^{\delta_1(t-mh)}-1))\\
+2(t-mh)\|P_{12}^{(m)}\|\gamma_{21}^{(m)}N_1M_2e^{(t-mh)(\delta_1+\mu_2)}\Big)\\
+\gamma_{21}^{(m)}\Big(\|P_{22}^{(m)}\|(\frac{\|A_{11}\|M_1N_2}{\mu_1}e^{\delta_2(t-mh)}(e^{\mu_1(t-mh)}-1)
+\frac{N_2\|A_{22}\|}{\delta_2}(e^{\delta_2(t-mh)}-1))\\
+2(t-mh)\|P_{12}^{(m)}\|\gamma_{12}^{(m)}N_2M_1e^{(t-mh)(\delta_2+\mu_1)}\Big)
:=\psi_{12}^{(m)}(t).
\endgathered
\end{equation*}

Further, we consider
\begin{equation*}
\gathered
e^{A_{22}^{\T}(t-mh)}(A_{12}^{\T}(t)P_{12}(t)+P_{21}(t)A_{12}(t))e^{A_{22}(t-mh)}
-(A_{12}^{\T}(t)P_{12}^{(m)}+P_{21}^{(m)}A_{12}(t))\\
=e^{A_{22}^{\T}(t-mh)}A_{12}^{\T}(t)P_{12}(t)e^{A_{22}(t-mh)}-
A_{12}^{\T}(t)P_{12}^{(m)}\\
+e^{A_{22}^{\T}(t-mh)}P_{21}(t)A_{12}(t)e^{A_{22}(t-mh)}-
P_{21}^{(m)}A_{12}(t)
\endgathered
\end{equation*}
Taking into account \eqref{5.6} we obtain
\begin{equation*}
\gathered
e^{A_{22}^{\T}(t-mh)}A_{12}^{\T}(t)P_{12}(t)e^{A_{22}(t-mh)}-
A_{12}^{\T}(t)P_{12}^{(m)}\\
=(e^{A_{22}^{\T}(t-mh)}A_{12}^{\T}(t)e^{-A_{11}^{\T}(t-mh)}-A_{12}^{\T}(t))P_{12}^{(m)}\\
-e^{A_{22}^{\T}(t-mh)}A_{12}^{\T}(t)e^{-A_{11}^{\T}(t-mh)}\int\limits_{mh}^t(P_{11}^{(m)}A_{12}(s)+A_{21}^{\T}(s)P_{22}^{(m)})\,ds
\endgathered
\end{equation*}
Consequently, applying \eqref{5.8***} we obtain
\begin{equation*}
\gathered
\|e^{A_{22}^{\T}(t-mh)}A_{12}^{\T}(t)P_{12}(t)e^{A_{22}(t-mh)}-
A_{12}^{\T}(t)P_{12}^{(m)}\|\\
\le\gamma_{12}^{(m)}\Big(\|P_{12}^{(m)}\|(\frac{\|A_{22}\|M_2N_1}{\mu_2}e^{\delta_1(t-mh)}(e^{\mu_2(t-mh)}-1)
+\frac{N_1\|A_{11}\|}{\delta_1}(e^{\delta_1(t-mh)}-1))\\
+(t-mh)(\|P_{22}^{(m)}\|\gamma_{21}^{(m)}+\gamma_{12}^{(m)}\|P_{11}^{(m)}\|)M_2N_1e^{(t-mh)(\mu_2+\delta_1)}\Big)
\endgathered
\end{equation*}
Thus, 
\begin{equation*}
\gathered
\|e^{A_{22}^{\T}(t-mh)}(A_{12}^{\T}(t)P_{12}(t)+P_{21}(t)A_{12}(t))e^{A_{22}(t-mh)}\\
-(A_{12}^{\T}(t)P_{12}^{(m)}+P_{21}^{(m)}A_{12}(t))\|\\
\le 2\gamma_{12}^{(m)}\Big(\|P_{12}^{(m)}\|(\frac{\|A_{22}\|M_2N_1}{\mu_2}e^{\delta_1(t-mh)}(e^{\mu_2(t-mh)}-1)
+\frac{N_1\|A_{11}\|}{\delta_1}(e^{\delta_1(t-mh)}-1))\\
+(t-mh)(\|P_{22}^{(m)}\|\gamma_{21}^{(m)}+\gamma_{12}^{(m)}\|P_{11}^{(m)}\|)M_2N_1e^{(t-mh)(\mu_2+\delta_1)}\Big):=\psi_{22}^{(m)}(t)
\endgathered
\end{equation*}
We define matrices
$\Psi_m(t)=(\psi_{ij}^{(m)}(t))_{i,j=1,2}$, then
\begin{equation}\label{DerLF}
\gathered
\dot v(t,x_1(t),x_2(t))\le \zeta_m^{\T}(t)\Psi_m(t)\zeta_m(t),
\endgathered
\end{equation}
where $\zeta_m(t)=(\|z_{1m}(t)\|,\|z_{2m}(t)\|)^{\T}$.
We represent $\Psi_m(t)$ in the following form
\begin{equation*}
\gathered
\Psi_m(t)=\gamma_{12}^{(m)}\Big(\|A_{22}\|M_2N_1e^{\delta_1(t-mh)}\frac{e^{\mu_2(t-mh)}-1}{\mu_2}+N_1\|A_{11}\|\frac{e^{\delta_1(t-mh)}-1}{\delta_1}\Big)\Upsilon_{1}\\
+\gamma_{21}^{(m)}\Big(\|A_{11}\|M_1N_2e^{\delta_2(t-mh)}\frac{e^{\mu_1(t-mh)}-1}{\mu_1}
+N_2\|A_{22}\|\frac{e^{\delta_2(t-mh)}-1}{\delta_2}\Big)\Upsilon_{2}\\
+2\gamma_{12}^{(m)}N_1M_2(t-mh)e^{(t-mh)(\mu_2+\delta_1)}\wt\Upsilon_{1}
+2\gamma_{21}^{(m)}N_2M_1(t-mh)e^{(t-mh)(\mu_1+\delta_2)}\wt\Upsilon_{2},
\endgathered
\end{equation*}
where
\begin{equation*}
\gathered
\Upsilon_{1}=
\begin{pmatrix}
0&\|P_{11}^{(m)}\|\\
\|P_{11}^{(m)}\|&2\|P_{12}^{(m)}\|
\end{pmatrix},\quad
\Upsilon_{2}=
\begin{pmatrix}
2\|P_{12}^{(m)}\|&\|P_{22}^{(m)}\|\\
\|P_{22}^{(m)}\|&0
\end{pmatrix}
\endgathered
\end{equation*}

\begin{equation*}
\gathered
\wt\Upsilon_{1}=
\begin{pmatrix}
0&2\|P_{12}^{(m)}\|\gamma_{21}^{(m)}\\
2\|P_{12}^{(m)}\|\gamma_{21}^{(m)}&(\|P_{11}^{(m)}\|\gamma_{12}^{(m)}+\|P_{22}^{(m)}\|\gamma_{21}^{(m)})
\end{pmatrix},\quad
\endgathered
\end{equation*}
\begin{equation*}
\gathered
\wt\Upsilon_{2}=
\begin{pmatrix}
(\|P_{11}^{(m)}\|\gamma_{12}^{(m)}+\|P_{22}^{(m)}\|\gamma_{21}^{(m)})&2\|P_{12}^{(m)}\|\gamma_{12}^{(m)}\\
2\|P_{12}^{(m)}\|\gamma_{12}^{(m)}&0
\end{pmatrix},\quad
\endgathered
\end{equation*}
Direct calculations show $\|\Upsilon_1\|=\eta_{11}^{(m)}$,
$\|\Upsilon_2\|=\eta_{22}^{(m)}$, $\|\wt\Upsilon_1\|=\eta_{12}^{(m)}$,
$\|\wt \Upsilon_2\|=\eta_{21}^{(m)}$.

For $\Psi_m(t)$ the following estimates hold
\begin{equation*}
\gathered
\|\Psi_m(t)\|\le \alpha_{11}^{(m)}e^{\delta_1(t-mh)}\frac{e^{\mu_2(t-mh)}-1}{\mu_2}
+\alpha_{12}^{(m)}\frac{e^{\delta_1(t-mh)}-1}{\delta_1}\\
+\alpha_{21}^{(m)}e^{\delta_2(t-mh)}\frac{e^{\mu_1(t-mh)}-1}{\mu_1}
+\alpha_{22}^{(m)}\frac{e^{\delta_2(t-mh)}-1}{\delta_2}
\\
+2(\gamma_{12}^{(m)}\eta_{12}^{(m)}N_1M_2(t-mh)e^{(t-mh)(\mu_2+\delta_1)}
+\gamma_{21}^{(m)}\eta_{21}^{(m)}N_2M_1(t-mh)e^{(t-mh)(\mu_1+\delta_2)}).
\endgathered
\end{equation*}
Therefore, for $\int_{mh}^{(m+1)h}\|\Psi_m(s)\|\,ds$ the following estimate holds
\begin{equation*}
\gathered
\int\limits_{mh}^{(m+1)h}\|\Psi_m(s)\|\,ds\le
\frac{\alpha_{11}^{(m)}}{\mu_2}\Big(\frac{e^{(\mu_2+\delta_1)h}-1}{\mu_2+\delta_1}-
\frac{e^{\delta_1h}-1}{\delta_1}\Big)
+\frac{\alpha_{12}^{(m)}}{\delta_1}\Big(\frac{e^{h\delta_1}-1}{\delta_1}-h\Big)\\
+\frac{\alpha_{21}^{(m)}}{\mu_1}\Big(\frac{e^{(\mu_1+\delta_2)h}-1}{\mu_1+\delta_2}-
\frac{e^{\delta_2h}-1}{\delta_2}\Big)
+\frac{\alpha_{22}^{(m)}}{\delta_2}\Big(\frac{e^{h\delta_2}-1}{\delta_2}-h\Big)+
\\
+2\Big(\gamma_{12}^{(m)}\eta_{12}^{(m)}N_1M_2\Big(\frac{he^{(\mu_2+\delta_1)h}}{\mu_2+\delta_1}-\frac{e^{(\mu_2+\delta_1)h}-1}{(\mu_2+\delta_1)^2}\Big)\\
+\gamma_{21}^{(m)}\eta_{21}^{(m)}N_2M_1\Big(\frac{he^{(\mu_1+\delta_2)h}}{\mu_1+\delta_2}-\frac{e^{(\mu_1+\delta_2)h}-1}{(\mu_1+\delta_2)^2}\Big)\Big)=\Theta_m(h).
\endgathered
\end{equation*}
From \eqref{DerLF}, taking into account the  Lemma 4.1, for
$t\in(mh,(m+1)h)$ we obtain 
\begin{equation*}
\gathered
\dot v(t,x_1(t),x_2(t))\le \zeta_m^{\T}\Psi_m(t)\zeta_m(t)
\le\|\Psi_m(t)\|(\|z_{1m}(t)\|^2+\|z_{2m}(t)\|^2)\\
\le
\frac{\|\Psi_m(t)\|}{\lambda_{\min}(\Pi_m)}v(t,x_1(t),x_2(t)),
\endgathered
\end{equation*}
Integrating this differential inequality we get
\begin{equation*}
\gathered
v((m+1)h,x_1((m+1)h),x_2((m+1)h))\\
\le e^{\int\limits_{mh}^{(m+1)h}\frac{\|\Psi_m(s)\|}{\lambda_{\min}(\Pi_m)}\,ds}
v(mh+0,x_1(mh+0),x_2(mh+0))\\
\le e^{\frac{\Theta_m(h)}{\lambda_{\min}(\Pi_m)}}
v(mh+0,x_1(mh+0),x_2(mh+0)),
\endgathered
\end{equation*}
which completes the proof of the Lemma.

\subsection{Conditions for the asymptotic stability}
We establish sufficient conditions for the asymptotic stability of \eqref{3.1} using  $\mathfrak V(t,x_1,x_2)$ constructed in the previous section and Lemma 4.2.
Recall that $P_0=(P_{ij}^{(0)})_{i,j=1,2}$, $P_N=(P_{ij}^{(N)})_{i,j=1,2}$ are block matrices, where $P_{ij}$ are defined by \eqref{5.1}-\eqref{5.3}.

{\bf Theorem 4.1.} Let $N\in\mathbb N$ and $P_0\succ 0$ be such that for $m=0,\dots,N-1$ the matrices $\Pi_m$ defined by \eqref{l4.1} are positive definite. If
\begin{equation}\label{teorem 4.1}
\gathered
Q:=\sum\limits_{m=0}^{N-1}\frac{\Theta_m(h)}{\lm_{\min}(\Pi_m)}+\ln\lm_{\max}(P_N^{-1}B^{\T}P_0B)<0,
\endgathered
\end{equation}
where $h=\frac{\theta}{N}$ and $\Theta_m(h)$ are defined in \eqref{Theta}, then \eqref{3.1} is asymptotically stable.

{\it Proof.} From  \eqref{5.9} it follows that for all $m=0,\dots,N-1$
\begin{equation*}
\gathered
v((m+1)h,x_1((m+1)h),x_2((m+1)h))\\
\le\exp\Big(\frac{\Theta_m(h)}{\lm_{\min}(\Pi_m)}\Big)v(mh+0,x_1(mh+0),x_2(mh+0)).
\endgathered
\end{equation*}
Therefore,
\begin{equation*}
\gathered
v(\theta,x_1(\theta),x_2(\theta))\le\exp\Big(\sum\limits_{m=0}^{N-1}\frac{\Theta_m(h)}{\lm_{\min}(\Pi_m)}\Big)v(0+0,x_1(0+0),x_2(0+0)).
\endgathered
\end{equation*}
For $t=\theta$,
\begin{equation*}
\gathered
v(\theta+0,x_1(\theta+0),x_2(\theta+0))= x^{\T}(\theta+0)P_0x(\theta+0)\\
=x^{\T}(\theta)B^{\T}P_0Bx(\theta)=
(P_N^{1/2}x(\theta))^{\T}P_N^{-1/2}B^{\T}P_0BP_N^{-1/2}P_N^{1/2}x(\theta)\\
\le\lm_{\max}(P_N^{-1/2}B^{\T}P_0BP_{N}^{-1/2})\|P_N^{1/2}x(\theta)\|^2=
\lm_{\max}(P_N^{-1}B^{\T}P_0B)v(\theta,x_1(\theta),x_2(\theta)).
\endgathered
\end{equation*}
We recall that here $x(\theta)=(x_1^{\T}(\theta),x_2^{\T}(\theta))^{\T}$.
Consequently,
\begin{equation*}
\gathered
v(\theta+0,x_1(\theta+0),x_2(\theta+0))\\
\le\lm_{\max}(P_N^{-1}B^{\T}P_0B)
\exp\Big(\sum\limits_{m=0}^{N-1}\frac{\Theta_m(h)}{\lm_{\min}(\Pi_m)}\Big)v(0+0,x_1(0+0),x_2(0+0))\\
=e^{Q}v(0+0,x_1(0+0),x_2(0+0)).
\endgathered
\end{equation*}
Due to the periodicity in $t$ of \eqref{3.1} and $\mathfrak V(t,.,.)$ for any $k\in\Bbb Z_+$ the following inequality holds
\begin{equation*}
\gathered
v((k+1)\theta+0,x_1((k+1)\theta+0),x_2((k+1)\theta+0))\\
\le
e^{Q}v(k\theta+0,x_1(k\theta+0),x_2(k\theta+0)).
\endgathered
\end{equation*}
Therefore,
\begin{equation*}
\gathered
v(k\theta+0,x_1(k\theta+0),x_2(k\theta+0))\le
e^{Qk}v(0+0,x_1(0+0),x_2(0+0)).
\endgathered
\end{equation*}
From Lemma 4.1 it follows that
\begin{equation*}
\gathered
\lm_{\min}(\Pi_0)\|x(k\theta+0)\|^2\le v(k\theta+0,x_1(k\theta+0),x_2(k\theta+0))\\
\le e^{Qk}v(0+0,x_1(0+0),x_2(0+0))\le e^{Qk}\lm_{\max}(\Xi_0)\|x(0+0)\|^2,
\endgathered
\end{equation*}
which is equivalent to
\begin{equation*}
\gathered
\|x(k\theta+0)\|\le
\sqrt{\frac{\lm_{\max}(\Xi_0)}{\lm_{\min}(\Pi_0)}}e^{Qk/2}\|x(0+0)\|.
\endgathered
\end{equation*}
By conditions of Theorem 4.1, $Q<0$, therefore, $\|x(k\theta+0)\|\to 0$ as $k\to\infty$ which proves the asymptotic stability of  \eqref{3.1}. The theorem is proved.

To check the condition \eqref{teorem 4.1}, it is necessary to calculate $\Theta_m(h)$. In addition, with the constants $\gamma_{ij}^{(m)}$, $m=0,\dots,N-1$ included in $\Theta_m(h)$, it is necessary to calculate the matrices $P_{ij}^{(m)}$, $i ,j=1,2$, $m=0,\dots,N$ using the recurrent formulas \eqref{5.1}-\eqref{5.3}. These formulas contain $e^{-A_{ii}h}$ and $\int_{mh}^{(m+1)h}A_{ij}(s)\,ds$. Such calculations are accessible to modern computing tools. It is intuitively clear that $h$ should be chosen such that $h\sup_{t\in[0,\theta], i,j=1,2}\|P_{ij}(t)\|\ll 1$ , i.e. so that the change of $P_{ij}(t)$ is rather small over the discretization period $(mh,(m+1)h]$.

\section{The case of a non-periodic system }
Consider the case of a linear impulsive system \eqref{3.1} when the dwell-times are subject to the conditions $\theta_1\le T_k\le\theta_2$, $\theta_1\ne\theta_2$. In this case, the system \eqref{3.1} is not periodic, since its solutions do not have the property of invariance with respect to the semigroup $\theta \Bbb Z_+$ and the Floquet theory is not applicable.


The Lyapunov function $v(t,x_1,x_2)$ has a quadratic form \eqref{5.7}, we will construct $P_{ij}(t)$, $i,j=1,2$  step by step on every interval $[\tau_k,\tau_{k+1}]$. If $\tau_k \notin h\Bbb{Z}$, then on the interval $(\tau_k,d_kh]$, where $d_kh$ is the grid node nearest on the left to $\tau_k$,  we choose $P_{ij}(t)$ as a constant on the interval $(\tau_k,d_kh]$. Between the nodes we construct $P_{ij}(t)$ similarly to the periodic case. Finally, if $\tau_{k+1} \notin h\Bbb{Z }$, then on the interval $(\varkappa_k h,\tau_{k+1}]$ ($\varkappa_kh$ is the grid node nearest on the right to $\tau_{k+1}$) we choose $P_{ij}(t)$ again as a constants, i.e. $P_{ij}(t)=P_{ij}(\varkappa_kh-0)$, for $t \in (\varkappa_kh,\tau_{k+1}]$. Thus, it becomes necessary to estimate the derivative of the Lyapunov function $v(t,x_1,x_2)$ on the time intervals $(\tau_k,d_kh]$, $(\varkappa_kh,\tau_{k+1}]$. For this we will introduce Assumption 5.1 and Lemma 5.1.
 
\subsection{Construction of the Lyapunov function} To construct a Lyapunov function, the discretization method is also used here. The discretization parameters are: the number of nodes $N\in\Bbb N$,  $N\ge 2$ and $h=\frac{\theta}{N}$. We denote $N_3=-\Big[\frac{2h-\theta_1}{h}\Big]$,
$N_4=\Big[\frac{\theta_2}{h}\Big]$. In this case, additional assumptions regarding the connections between the subsystems are required.

{\bf Assumption 5.1.} There are positive constants $l_{ij}^{(m)}$, $i,j=1,2$, $i\ne j$, $m=0,\dots,N-1$ such that
\begin{equation*}
\gathered
\sup\limits_{t\in(mh,(m+1)h]}\|A_{ij}(t)-A_{ij}(mh)\|\le l_{ij}^{(m)}h.
\endgathered
\end{equation*}
For any $m\in\Bbb Z$, let $l_{ij}^{(m)}:=l_{ij}^{(\varrho)}$, where $\varrho$ is the remainder of dividing $m$ by $N$.
Let $\mathfrak{l}_m:=\sqrt{(l_{12}^{(m)})^2+(l_{21}^{(m)})^2}$.

We denote the constant matrices
\begin{equation*}
\gathered
A_m=
\begin{pmatrix}
A_{11}&A_{12}(mh)\\
A_{21}(mh)&A_{22}
\end{pmatrix}, \quad m\in\Bbb Z.
\endgathered
\end{equation*}

Let $P_0=(P_{ij}^{(0)})_{i,j=1,2}$, $P_{ij}^{(0)}=(P_{ij}^{(0)})^{\T}$ be a positive definite symmetric block matrix.
We define sequences of block matrices $P_{m}^{(l)}=(P_{ij}^{(m,l)})_{i,j=1,2}$, $l=0,\dots,N-1$, $m=0,\dots,N_4-1$ as follows $P_{ij}^{(0,l)}\equiv P_{ij}^{(0)} $
\begin{equation}\label{6.1}
\gathered
P_{11}^{(m+1,l)}= e^{-A_{11}^{\T}h}(P_{11}^{(m,l)}\\
-\int\limits_{(m+l)h}^{(m+l+1)h}(P_{12}^{(m,l)}A_{21}(s)+A_{21}^{\T}(s)P_{21}^{(m,l)})\,ds)e^{-A_{11}h},\\
\endgathered
\end{equation}
\begin{equation}\label{6.2}
\gathered
P_{22}^{(m+1,l)}=e^{-A_{22}^{\T}h}(P_{22}^{(m,l)}\\
-\int\limits_{(m+l)h}^{(m+l+1)h}(A_{12}^{\T}(s)P_{12}^{(m,l)}+P_{21}^{(m,l)}A_{12}(s))\,ds)e^{-A_{22}h},\\\\
\endgathered
\end{equation}
\begin{equation}\label{6.3}
\gathered
P_{12}^{(m+1,l)}= e^{-A_{11}^{\T}h}(P_{12}^{(m,l)}\\
-\int\limits_{(m+l)h}^{(m+l+1)h}(P_{11}^{(m,l)}A_{12}(s)+A_{21}^{\T}(s)P_{22}^{(m,l)})\,ds)e^{-A_{22}h}.
\endgathered
\end{equation}
Next, we define the matrices $P_{ij}(t)$, $i,j=1,2$ sequentially on the intervals $(\tau_k,\tau_{k+1}]$.
Let $\wt\tau_k=\tau_k-\Big[\frac{\tau_k}{\theta}\Big]\theta\in[0,\theta)$,
$l_k:=\Big[\frac{\wt\tau_k}{h}\Big]+1$, $d_k:=\Big[\frac{\tau_k}{h}\Big]+1$, $\varkappa_k=\Big[\frac{\tau_{k+1}}{h}\Big]$. It is easy to show that
\begin{equation*}
\gathered
(\tau_k,\tau_{k+1}]=(\tau_k,d_kh]\cup\bigcup\limits_{l=d_k}^{\varkappa_k-1}(lh,(l+1)h]\cup(\varkappa_kh,\tau_{k+1}].
\endgathered
\end{equation*}
Let $P_{ij}(t)=P_{ij}^{(0)}$ for $t\in(\tau_k,d_kh]$. On each interval $((m+d_k)h,(m+1+d_k)h]$, $m=0,\dots,\varkappa_k-d_k-1$, we put
\begin{equation}\label{6.4}
\gathered
P_{11}(t)= e^{-A_{11}^{\T}(t-(m+d_k)h)}(P_{11}^{(m,l_k)}\\
-\int\limits_{(m+d_k)h}^t(P_{12}^{(m,l_k)}A_{21}(s)+A_{21}^{\T}(s)P_{21}^{(m,l_k)})\,ds)e^{-A_{11}(t-(m+d_k)h)},
\endgathered
\end{equation}
\begin{equation}\label{6.5}
\gathered
P_{22}(t)=e^{-A_{22}^{\T}(t-(m+d_k)h)}(P_{22}^{(m,l_k)}\\
-\int\limits_{(m+d_k)h}^t(A_{12}^{\T}(s)P_{12}^{(m,l_k)}+P_{21}^{(m,l_k)}A_{12}(s))\,ds)e^{-A_{22}(t-(m+d_k)h)},
\endgathered
\end{equation}
\begin{equation}\label{6.6}
\gathered
P_{12}(t)= e^{-A_{11}^{\T}(t-(m+d_k)h)}(P_{12}^{(m,l_k)}\\
-\int\limits_{(m+d_k)h}^t(P_{11}^{(m,l_k)}A_{12}(s)+A_{21}^{\T}(s)P_{22}^{(m,l_k)})\,ds)e^{-A_{22}(t-(m+d_k)h)}.
\endgathered
\end{equation}
For $t\in(\varkappa_kh,\tau_{k+1}]$, let $P_{ij}(t)=P_{ij}^{(\varkappa_k-d_k,l_k)}$.

{\it Remark.} Note that $P(t)$ is continuous on $(\tau_k,\tau_{k+1})$ and left-continuous at $t=\tau_{k+1}$. Indeed, due to the $\theta$-periodicity of
$A_{ij}(t)$, the matrices $P_{ij}(t)$ satisfy 
\begin{equation}\label{6.7'}
\gathered
P_{ij}((m+d_k)h+0)=P_{ij}^{(m,l_k)}, \quad m=0,\dots,\varkappa_k-d_k-1,
\endgathered
\end{equation}
\begin{equation}\label{6.8'}
\gathered
P_{ij}((m+d_k)h-0)=P_{ij}^{(m,l_k)}, \quad m=0,\dots,\varkappa_k-d_k-1.
\endgathered
\end{equation}

Let us prove \eqref{6.7'} and \eqref{6.8'}. We restrict ourselves to the case $(i,j)=(1,1)$, since other cases one can consider in a similar way. In case $m=0$ the formulas \eqref{6.7'} and \eqref{6.8'} are obvious. From the equality \eqref{6.4} it is easy to show that \eqref{6.7'} is true for all $m = 1, \dots, \varkappa_k-d_k-1$.

Next we prove \eqref{6.8'}. From \eqref{6.4} we have
\begin{equation}\label{6.9'}
\gathered
P_{11}((m+d_k)h-0)= e^{-A_{11}^{\T}h}(P_{11}^{(m-1,l_k)}\\
-\int\limits_{(m+d_k-1)h}^{(m+d_k)h}(P_{12}^{(m-1,l_k)}A_{21}(s)+A_{21}^{\T}(s)P_{21}^{(m-1,l_k)})\,ds)e^{-A_{11}h},
\endgathered
\end{equation}
In the integral \eqref{6.9'} make the change of variables $\widetilde{s}:=s+(l_k-d_k)h$, then
\begin{equation}\label{6.10'}
\gathered
P_{11}((m+d_k)h-0)= e^{-A_{11}^{\T}h}(P_{11}^{(m-1,l_k)}\\
-\int\limits_{(m+l_k-1)h}^{(m+l_k)h}(P_{12}^{(m-1,l_k)}A_{21}(\widetilde{s}+(l_k-d_k)h)\\
+A_{21}^{\T}(\widetilde{s}+(l_k-d_k)h)P_{21}^{(m-1,l_k)})\,ds)e^{-A_{11}h},
\endgathered
\end{equation}

Note that $(d_k-l_k)h \in \theta\Bbb{Z}$. Indeed, the definition of the function $x\mapsto [x]$ implies the inequalities
\begin{equation}\label{6.11'}
\gathered
\Big [\frac{\tau_k}{h}\Big]\le\frac{\tau_k}{h}<\Big [\frac{\tau_k}{h}\Big]+1, \quad \Big [\frac{\widetilde{\tau}_k}{h}\Big]\le\frac{\widetilde{\tau}_k}{h}<\Big [\frac{\widetilde{\tau}_k}{h}\Big]+1.
\endgathered
\end{equation}

Therefore, from \eqref{6.11'} it follows that
\begin{equation*}
\gathered
\frac{\tau_k-\widetilde{\tau}_k}{h}-1<\Big[\frac{\tau_k}{h} \Big]-\Big[\frac{\widetilde{\tau}_k}{h} \Big]<\frac{\tau_k-\widetilde{\tau}_k}{h}+1.
\endgathered
\end{equation*}

From the definition of $\widetilde{\tau}_k$ follows $\frac{\tau_k- \widetilde{\tau}_k}{h}=[\frac{\tau_k}{\theta}] \frac{\theta }{h}=[\frac {\tau_k}{\theta}]N \in \Bbb{Z}$, hence
\begin{equation}\label{6.12'}
\gathered
\Big[\frac{\tau_k}{\theta} \Big]N-1<\Big[\frac{\tau_k}{h} \Big]-\Big[\frac{\widetilde{\tau}_k}{h} \Big]<\Big[\frac{\widetilde{\tau}_k}{\theta} \Big]N+1.
\endgathered
\end{equation}

Since $[\frac{\tau_k}{h}]-[\frac{\widetilde{\tau}_k}{h}]\in \Bbb{Z}$, then $[\frac{\tau_k}{h}]-[\frac{\widetilde{\tau}_k}{h}]=[\frac{\tau_k}{\theta}]N$ and $h([\frac{\tau_k}{h}]-[\frac{\widetilde{\tau}_k}{h}])=[\frac{\tau_k}{\theta}]Nh=[\frac{\tau_k}{\theta}]\theta \in \theta\Bbb{Z}$, which means $(d_k-l_k)h \in \theta \Bbb{Z}$. In turn, $A_{12}(s)$ and $A_{21}(s)$ are periodic functions, therefore $A_{12}(s+(d_k-l_k)h)=A_{12}(s)$ and $A_{21}(s+(d_k-l_k)h)=A_{21}(s)$, then for equality \eqref{6.10'}, taking into account \eqref{6.1}, $P_{11}((m+d_k)h-0)=P_{11}^{(m,l)}$.

We choose the elements $\mathfrak V(t,x)$ in the form $v_{ij}(t,x_i,x_j)=x_i^{\T}P_{ij}(t)x_j$, $i,j=1,2$.


The following lemma uses Assumption 5.1 and is needed to estimate the derivative of the Lyapunov function $v(t,x_1,x_2)$ on $[\tau_k, d_kh)$ and $(\varkappa_k h, \tau_{k+1}]$, where $P_{ij}(t)$ is constant.

{\bf Lemma 5.1.} Let $\mathfrak{l}_m:=\sqrt{(l_{12}^{(m)})^2+(l_{21}^{(m)})^2}$, then
\begin{equation*}
\gathered
\dot v(t,x(t))\le \lambda_{\max}(P_0^{-1}(A_{l_k-1}^{\T}P_0+P_0A_{l_k-1}+2h\mathfrak{l}_{l_k-1}\|P_0\|I_n))v(t,x(t)) ,\quad t\in(\tau_k,d_kh],
\endgathered
\end{equation*}
\begin{equation*}
\gathered
\dot v(t,x(t))\le \lambda_{\max}((P_{\varkappa_k-d_k}^{(l_k)})^{-1}(A_{l_k+\varkappa_k-d_k}^{\T}P_{\varkappa_k-d_k}^{(l_k)}+P_{\varkappa_k-d_k}^{(l_k)}
A_{l_k+\varkappa_k-d_k}\\
+2h\mathfrak{l}_{l_k+\varkappa_k-d_k}\|P_{\varkappa_k-d_k}^{(l_k)}\|I_n))v(t,x(t)), \quad t\in(\varkappa_kh,\tau_{k+1}].
\endgathered
\end{equation*}

{\it Proof.} Let $t\in(\tau_k,d_kh]$, then, taking into account that $(l_k-d_k)h \in \theta \Bbb{Z}$
\begin{equation}\label{IIder}
\gathered
\dot v(t,x(t))=x^{\T}(t)(A^{\T}(t)P(t)+P(t)A(t))x(t)
=x^{\T}(t)(A^{\T}(t)P_0+P_0A(t))x(t)\\=
x^{\T}(t)(A^{\T}((d_k-1)h)P_0+P_0A((d_k-1)h))x(t)\\
+x^{\T}(t)((A(t)-A((d_k-1)h))^{\T}P_0+P_0(A(t)-A((d_k-1)h)))x(t)\\
\le x^{\T}(t)(A^{\T}((l_k-1)h)P_0+P_0A((l_k-1)h))x(t)+2h\mathfrak{l}_{l_k-1}\|P_0\| x^{\T}(t)x(t)\\
\le(P_0^{1/2}x(t))^{\T}P_0^{-1/2}(A^{\T}((l_k-1)h)P_0+P_0A((l_k-1)h)
+2h\mathfrak{l}_{l_k-1}\|P_0\|I_n)P_0^{-1/2}(P_0^{1/2}x(t))\\
\le\lambda_{\max}(P_0^{-1}(A_{l_k-1}^{\T}P_0+P_0A_{l_k-1}+2h\mathfrak{l}_{l_k-1}\|P_0\|I_n))v(t,x(t)).
\endgathered
\end{equation}
The second part of the lemma can be proved similarly. The lemma is proved.


Lemma 5.2 is similar to Lemma 4.1 and reduces the verification of the positive definiteness of the Lyapunov function to a finite number of inequalities.

{\bf Lemma 5.2.} If $z_{im}=e^{-A_{ii}(t-mh)}x_i$ for $t\in(mh,(m+1)h],$ $i=1,2$, then
\begin{equation}\label{5.8}
\gathered
\lm_{\min}(\Pi_m^{(l)})\|z_m\|^2\le v(t,x_1,x_2)\le\lm_{\max}(\Xi_m^{(l)})\|z_m\|^2,\\
\text{for all} \quad t\in(mh,(m+1)h],\quad m=0,\dots,N-1,
\endgathered
\end{equation}
where $z_m=(z_{1m}^{\T},z_{2m}^{\T})^{\T}$, $\|z_m\|^2=\|z_{1m}\|^2+\|z_{2m}\|^2$, $\Pi_m^{(l)}=(\pi_{ij}^{(m,l)})_{i,j=1,2}$, $\Xi_m^{(l)}=(\xi_{ij}^{(m,l)})_{i,j=1,2}$ are block matrices with the elements
\begin{equation}\label{33*}
\gathered
\pi_{11}^{(m,l)}=P_{11}^{(m,l)}-h(2\gamma_{21}^{(m+l)}\|P_{12}^{(m,l)}\|+(\gamma_{12}^{(m+l)}\|P_{11}^{(m,l)}\|+\gamma_{21}^{(m+l)}\|P_{22}^{(m,l)}\|))I_{n_1},\\
\pi_{22}^{(m)}=P_{22}^{(m,l)}-h(2\gamma_{12}^{(m+l)}\|P_{12}^{(m,l)}\|+(\gamma_{12}^{(m+l)}\|P_{11}^{(m,l)}\|+\gamma_{21}^{(m+l)}\|P_{22}^{(m,l)}\|))I_{n_2},\\
\pi_{12}^{(m)}=P_{12}^{(m,l)},\quad \pi_{21}^{(m,l)}=P_{21}^{(m,l)}
\endgathered
\end{equation}
\begin{equation*}
\gathered
\xi_{11}^{(m,l)}=P_{11}^{(m,l)}+h(2\gamma_{21}^{(m+l)}\|P_{12}^{(m,l)}\|+(\gamma_{12}^{(m+l)}\|P_{11}^{(m,l)}\|+\gamma_{21}^{(m+l)}\|P_{22}^{(m,l)}\|))I_{n_1},\\
\xi_{22}^{(m,l)}=P_{22}^{(m,l)}+h(2\gamma_{12}^{(m+l)}\|P_{12}^{(m,l)}\|+(\gamma_{12}^{(m+l)}\|P_{11}^{(m,l)}\|+\gamma_{21}^{(m+l)}\|P_{22}^{(m,l)}\|))I_{n_2},\\
\xi_{12}^{(m,l)}=P_{12}^{(m,l)},\quad \xi_{21}^{(m,l)}=P_{21}^{(m,l)}.
\endgathered
\end{equation*}
The proof is similar to the proof of Lemma 4.1.

We denote
\begin{equation*}
\gathered
\eta_{11}^{(m,l)}:=\sqrt{\|P_{11}^{(m,l)}\|^2+\|P_{12}^{(m,l)}\|^2}+\|P_{12}^{(m,l)}\|,\\
\eta_{22}^{(m,l)}:=\sqrt{\|P_{22}^{(m,l)}\|^2+\|P_{12}^{(m,l)}\|^2}+\|P_{12}^{(m,l)}\|
\endgathered
\end{equation*}
\begin{equation*}
\gathered
\eta_{12}^{(m,l)}:=\frac{1}{2}(\|P_{11}^{(m,l)}\|\gamma_{12}^{(m+l)}+\|P_{22}^{(m,l)}\|\gamma_{21}^{(m+l)}\\
+\sqrt{(\|P_{11}^{(m,l)}\|\gamma_{12}^{(m+l)}+\|P_{22}^{(m,l)}\|\gamma_{21}^{(m+l)})^2
+16(\gamma_{21}^{(m,l)})^2\|P_{12}^{(m,l)}\|^2}).
\endgathered
\end{equation*}
\begin{equation*}
\gathered
\eta_{21}^{(m,l)}:=\frac{1}{2}(\|P_{11}^{(m,l)}\|\gamma_{12}^{(m+l)}+\|P_{22}^{(m,l)}\|\gamma_{21}^{(m+l)}\\
+\sqrt{(\|P_{11}^{(m,l)}\|\gamma_{12}^{(m+l)}+\|P_{22}^{(m,l)}\|\gamma_{21}^{(m+l)})^2
+16(\gamma_{12}^{(m+l)})^2\|P_{12}^{(m,l)}\|^2}).
\endgathered
\end{equation*}
\begin{equation*}
\gathered
\alpha_{11}^{(m,l)}:=
\gamma_{12}^{(m+l)}\|A_{22}\|N_1M_2\eta_{11}^{(m,l)},\quad
\alpha_{12}^{(m,l)}:=\gamma_{12}^{(m+l)}\|A_{11}\|N_1\eta_{11}^{(m,l)},
\endgathered
\end{equation*}
\begin{equation*}
\gathered
\alpha_{21}^{(m,l)}:=
\gamma_{21}^{(m+l)}\|A_{11}\|N_2M_1\eta_{22}^{(m,l)},\quad
\alpha_{22}^{(m,l)}:=\gamma_{21}^{(m+l)}\|A_{22}\|N_2\eta_{22}^{(m,l)},
\endgathered
\end{equation*}
Let
\begin{equation}\label{Theta*}
\gathered
\Theta_{m}^{(l)}(h):=\frac{\alpha_{11}^{(m,l)}}{\mu_2}\Big(\frac{e^{(\mu_2+\delta_1)h}-1}{\mu_2+\delta_1}-
\frac{e^{\delta_1h}-1}{\delta_1}\Big)
+\frac{\alpha_{12}^{(m,l)}}{\delta_1}\Big(\frac{e^{h\delta_1}-1}{\delta_1}-h\Big)\\
+\frac{\alpha_{21}^{(m,l)}}{\mu_1}\Big(\frac{e^{(\mu_1+\delta_2)h}-1}{\mu_1+\delta_2}-
\frac{e^{\delta_2h}-1}{\delta_2}\Big)
+\frac{\alpha_{22}^{(m,l)}}{\delta_2}\Big(\frac{e^{h\delta_2}-1}{\delta_2}-h\Big)+
\\
+2\Big(\gamma_{12}^{(m+l)}\eta_{12}^{(m,l)}N_1M_2\Big(\frac{he^{(\mu_2+\delta_1)h}}{\mu_2+\delta_1}-\frac{e^{(\mu_2+\delta_1)h}-1}{(\mu_2+\delta_1)^2}\Big)\\
+\gamma_{21}^{(m+l)}\eta_{21}^{(m,l)}N_2M_1\Big(\frac{he^{(\mu_1+\delta_2)h}}{\mu_1+\delta_2}-\frac{e^{(\mu_1+\delta_2)h}-1}{(\mu_1+\delta_2)^2}\Big)\Big).
\endgathered
\end{equation}


Lemma 5.3 is similar to Lemma 4.2 and allows us to estimate the variations of the Lyapunov function $v(t,x_1,x_2)$ on each discritization interval $((m+d_k)h, (m+d_k+1)h)$.

{\bf Lemma 5.3}
Let $\Pi_m^{(l)}$ be positive definite for $m=0,\dots,\varkappa_k-d_k-1$, $l=0,\dots,N-1$, then for $t\in((m+d_k)h,(m+d_k+1)h)$, we have

\begin{equation}\label{6.9}
\gathered
v((m+d_k+1)h,x_1((m+d_k+1)h),x_2((m+d_k+1)h))\\
\le e^{\frac{\Theta_m^{(l)}(h)}{\lm_{\min}(\Pi_m^{(l)})}}v((m+d_k)h+0,x_1((m+d_k)h+0),x_2((m+d_k)h+0))
\endgathered
\end{equation}

The proof of Lemma 5.3 is similar to the proof of Lemma 4.2.

\subsection{Conditions for the asymptotic stability}

The Lyapunov function $v(t,x_1,x_2)$  constructed in the previous subsection and and Lemmas 5.1-5.3 allows us to formulate sufficient conditions for the asymptotic stability of \eqref{3.1} when $\theta_1 \ne \theta_2$.

{\bf Theorem 5.1.} Let $N\in\mathbb N$ and $P_0\succ 0$ be such that for
$m=0,\dots,N_4$, $l=0,\dots,N-1$, the matrices $\Pi_m^{(l)}$ defined by \eqref{33*} are positive definite and for all $M\in\mathbb N$, such that $N_3\le M\le N_4-1$ the following inequality hold
\begin{equation}\label{theorem5.1}
\gathered
h\lambda_{\max}^+(P_0^{-1}(A_{l-1}^{\T}P_0+P_0A_{l-1}+2h\mathfrak{l}_{l-1}\|P_0\|I_n))\\
+h\lambda_{\max}^+((P_{M}^{(l)})^{-1}(A_{l+M}^{\T}P_{M}^{(l)}+P_{M}^{(l)}
A_{l+M}+2h\mathfrak{l}_{l+M}\|P_{M}^{(l)}\|I_n))\\
+\sum\limits_{m=0}^{M-1}\frac{\Theta_m^{(l)}(h)}{\lm_{\min}(\Pi_m^{(l)})}+\ln\lm_{\max}((P_M^{(l)})^{-1}B^{\T}P_0B)<0,
\endgathered
\end{equation}
where $h=\frac{\theta}{N}$, $P_{ij}^{(m,l)}$ are defined in \eqref{6.1}-\eqref{6.3} and $\Theta_m^{(l)}(h)$ in \eqref{Theta*}. Then \eqref{3.1} is asymptotically stable, if $T_k \in [\theta_1,\theta_2]$. 

{\it Proof.} Consider the variation of the Lyapunov function on the interval $(\tau_k,\tau_{k+1}]$. As a consequence of Lemma 5.1, we have
\begin{equation}\label{6.10}
\gathered
v(hd_k,x(hd_k))\\
\le\exp\Big(\lambda_{\max}(P_0^{-1}(A_{l_k-1}^{\T}P_0+P_0A_{l_k-1}+2h\mathfrak{l}_{l_k-1}\|P_0\|I_n))(hd_k-\tau_k)\Big)
 v(\tau_k+0,x(\tau_k+0))\\
\le\exp\Big(\lambda_{\max}^+(P_0^{-1}(A_{l_k-1}^{\T}P_0+P_0A_{l_k-1}+2h\mathfrak{l}_{l_k-1}\|P_0\|I_n))h\Big)v(\tau_k+0,x(\tau_k+0)),
\endgathered
\end{equation}
\begin{equation}\label{6.11}
\gathered
v(\tau_{k+1},x(\tau_{k+1}))\\
\le\exp\Big(\lambda_{\max}((P_{\varkappa_k-d_k}^{(l_k)})^{-1}(A_{l_k+\varkappa_k-d_k}^{\T}P_{\varkappa_k-d_k}^{(l_k)}+P_{\varkappa_k-d_k}^{(l_k)}
A_{l_k+\varkappa_k-d_k}\\
+2h\mathfrak{l}_{l_k+\varkappa_k-d_k}\|P_{\varkappa_k-d_k}^{(l_k)}\|I_n))(\tau_{k+1}-h\varkappa_k)\Big)v(\tau_k+0,x(\tau_k+0))\\
\le\exp\Big(\lambda_{\max}^+((P_{\varkappa_k-d_k}^{(l_k)})^{-1}(A_{l_k+\varkappa_k-d_k}^{\T}P_{\varkappa_k-d_k}^{(l_k)}+P_{\varkappa_k-d_k}^{(l_k)}
A_{l_k+\varkappa_k-d_k}\\
+2h\mathfrak{l}_{l_k+\varkappa_k-d_k}\|P_{\varkappa_k-d_k}^{(l_k)}\|I_n))h\Big)v(\varkappa_kh+0,x(\varkappa_kh+0)).
\endgathered
\end{equation}
Lemma 5.3 implies the inequalities
\begin{equation}\label{6.12}
\gathered
v((m+1+d_k)h,x((m+1+d_k)h))\\
\le\exp\Big(\frac{\Theta_m^{(l)}(h)}{\lm_{\min}(\Pi_{m}^{(l_k)})}\Big)v((m+d_k)h+0,x((m+d_k)h+0))
\endgathered
\end{equation}
for all $m=0,\dots,\varkappa_k-d_k-1$.

At the moment of impulse action $t=\tau_{k+1}$, we get
\begin{equation}\label{6.13}
\gathered
v(\tau_{k+1}+0,x(\tau_{k+1}+0))\\=x^{\T}(\tau_{k+1}+0)P_0x(\tau_{k+1}+0)
=x^{\T}(\tau_{k+1})B^{\T}P_0Bx(\tau_{k+1})\\
=((P_{\varkappa_k-d_k}^{(l_k)})^{1/2}x(\tau_{k+1}))^{\T}(P_{\varkappa_k-d_k}^{(l_k)})^{-1/2}B^{\T}
P_0B(P_{\varkappa_k-d_k}^{(l_k)})^{-1/2}(P_{\varkappa_k-d_k}^{(l_k)})^{1/2}x(\tau_{k+1})\\
\le\lm_{\max}((P_{\varkappa_k-d_k}^{(l_k)})^{-1}B^{\T}P_0B)v(\tau_{k+1},x(\tau_{k+1}))
\endgathered
\end{equation}
Comparing the inequalities \eqref{6.10}---\eqref{6.13}, we obtain
\begin{equation}\label{6.14}
\gathered
v(\tau_{k+1}+0,x(\tau_{k+1}+0))\le e^{Q_{\varkappa_k-d_k}^{(l_k)}}v(\tau_{k}+0,x(\tau_{k}+0)),
\endgathered
\end{equation}
where
\begin{equation*}
\gathered
Q_{s}^{(l)}=h\lambda_{\max}^+(P_0^{-1}(A_{l-1}^{\T}P_0+P_0A_{l-1}+2h\mathfrak{l}_{l-1}\|P_0\|I_n))\\
+h\lambda_{\max}^+((P_{s}^{(l)})^{-1}(A_{l+s}^{\T}P_{s}^{(l)}+P_{s}^{(l)}
A_{l+s}+2h\mathfrak{l}_{l+s}\|P_{s}^{(l)}\|I_n))\\
+\sum\limits_{m=0}^{s-1}\frac{\Theta_m^{(l)}(h)}{\lm_{\min}(\Pi_m^{(l)})}+\ln\lm_{\max}((P_s^{(l)})^{-1}B^{\T}P_0B)
\endgathered
\end{equation*}
Note that from the condition $\theta_1\leq \theta_2$ it follows that $N_3\leq N_4-1$. By assumption about dwell-times, $\theta_1\le T_k\le\theta_2$. Therefore,
\begin{equation*}
\gathered
h(\varkappa_k-d_k)< T_k\le\theta_2,
\endgathered
\end{equation*}
and since $\varkappa_k-d_k\in\Bbb Z_+$, we get $\varkappa_k-d_k\le N_4-1$. On the other hand,
\begin{equation*}
\gathered
\theta_1\le T_k=\tau_{k+1}-\varkappa_kh+h(\varkappa_k-d_k)+d_kh-\tau_k\le h(\varkappa_k-d_k)+2h
\endgathered
\end{equation*}
Consequently, $\varkappa_k-d_k\ge\frac{\theta_1}{h}-2$, and since $\varkappa_k-d_k\in\Bbb Z_+$, we have
$\varkappa_k-d_k\ge N_3$. From the conditions of the theorem, it follows that
\begin{equation*}
\gathered
Q_{\max}=\max\limits_{l=0,\dots,N-1,N_3\le s\le N_4-1}Q_{s}^{(l)}<0.
\endgathered
\end{equation*}
It follows from the inequality \eqref{6.14} that for all $k\in\Bbb Z_+$,
\begin{equation}\label{6.15}
\gathered
v(\tau_{k+1}+0,x_1(\tau_{k+1}+0),x_2(\tau_{k+1}+0))\le e^{Q_{\max}}v(\tau_{k}+0,x_1(\tau_{k}+0),x_2(\tau_{k}+0)).
\endgathered
\end{equation}
From the inequalities \eqref{6.15}, we find
\begin{equation}\label{6.16}
\gathered
v(\tau_{k}+0,x(\tau_{k}+0))\le e^{kQ_{\max}}v(\tau_{0}+0,x(\tau_{0}+0)).
\endgathered
\end{equation}
Hence,
\begin{equation*}
\gathered
\lm_{\min}(P_0)\|x(\tau_{k}+0)\|^2\le v(\tau_{k}+0,x(\tau_{k}+0))\\
\le e^{kQ_{\max}}v(\tau_{0}+0,x(\tau_{0}+0))\le
e^{kQ_{\max}}\lm_{\max}(P_0)\|x(\tau_{0}+0)\|^2.
\endgathered
\end{equation*}
and
\begin{equation*}
\gathered
\|x(\tau_{k}+0)\|\le
e^{kQ_{\max}/2}\sqrt{\frac{\lm_{\max}(P_0)}{\lm_{\min}(P_0)}}\|x(\tau_{0}+0)\|.
\endgathered
\end{equation*}
Without loss of generality, we assume that $t_0<\tau_0$. Using the Gronwall-Bellman inequality, one can show that there exists a positive constant $C_{t_0}$ such that $\|x(\tau_{0}+0)\|\le C_{t_0}\|x_0\|$. It follows from the inequality $T_k<\theta_2$ that there exists a positive constant $C_2$ such that for all $t\in(\tau_k,\tau_{k+1}]$, the inequality $\|x(t)\|\le C_2\|x(\tau_{k}+0)\|$ holds.
Let $t\in(\tau_k,\tau_{k+1}]$. Then, the inequality
\begin{equation*}
\gathered
t-t_0=t-\tau_k+\sum\limits_{s=1}^kT_s+\tau_0-t_0\le k\theta_2+\theta_2+\tau_0-t_0
\endgathered
\end{equation*}
implies that
\begin{equation*}
\gathered
kQ_{\max}\le\frac{Q_{\max}}{\theta_2}(t-t_0)-\frac{(\theta_2+\tau_0-t_0)Q_{\max}}{\theta_2}.
\endgathered
\end{equation*}
Thus,
\begin{equation}\label{6.17}
\gathered
\|x(t)\|\le Ce^{-\beta(t-t_0)}\|x_0\|,\quad t\ge t_0,
\endgathered
\end{equation}
where
\begin{equation*}
\gathered
C=C_{t_0}C_2\exp\Big(-\frac{(\theta_2+\tau_0-t_0)Q_{\max}}{2\theta_2}\Big)\sqrt{\frac{\lm_{\max}(P_0)}{\lm_{\min}(P_0)}},\quad
\beta=-\frac{Q_{\max}}{2\theta_2}>0.
\endgathered
\end{equation*}
The estimate \eqref{6.17} implies the exponential stability of \eqref{3.1}. The theorem is proved.


To verify the conditions of Theorem 5.1, it is necessary to calculate $\Theta_m^{(l)}(h)$. In addition, with the constants $\gamma_{ij}^{(m+l)}$, $m=0,\dots,N-1$, $l=0,\dots,N-1$ included in $\Theta_m^{(l)}(h)$, it is necessary to calculate the matrices $P_{ij}^{(m,l)}$, $i ,j=1,2$, using the recurrent formulas \eqref{5.1}-\eqref{5.3} and check a finite number $N_4-N_3$ of inequalities \eqref{theorem5.1}.

\section{Comparison of results}

Here we compare the results obtained in the Section 4 with the known stability conditions for coupled time-variant systems. We restrict ourselves to a special case of high-frequency periodic functions, assuming dwell-time $T_k=\theta$ are constant and $\theta$ is a sufficiently small parameter. The value of this parameter we obtain from Theorem 4.1 applied to the case $N=1$. We will compare our results with small-gain conditions obtained on the basis of the ISS theory and the Lyapunov vector function. We will also consider the case when one of the independent subsystems is not stable, and the known approaches from the theory of stability of coupled systems are not applicable.

Consider  \eqref{3.1}
with $\widehat{A}_{ij}=\frac{1}{\theta}\int_{0}^{\theta}A_{ij}(t)\,dt$ for $i\ne j$.
Let $B=(B_{ij})_{i,j=1,2}$,
$
\widehat A=
\begin{pmatrix}
0&\widehat A_{12}\\
\widehat A_{21}&0
\end{pmatrix}
$
are block matrices. Consider the system of linear matrix inequalities
\begin{equation}\label{6.0}
\gathered
\diag\{e^{\theta A_{11}^{\T}},e^{\theta A_{22}^{\T}}\}B^{\T}P_0
B\diag\{e^{\theta A_{11}},e^{\theta A_{22}}\}\prec P_0-\theta(\widehat A^{\T}P_0+P_0\widehat A).
\endgathered
\end{equation}
Suppose it has a solution $P_0=(P_{ij}^{(0)})_{i,j=1,2}$ as a symmetric positive-definite matrix,
i.e. $P_{ij}^{(0)}=(P_{ji}^{(0)})^{\T}$. Let us define matrices
\begin{equation*}
\gathered
\pi_{11}^{(0)}=P_{11}^{(0)}-\theta(2\gamma_{21}^{(0)}\|P_{12}^{(0)}\|+(\gamma_{12}^{(0)}\|P_{11}^{(0)}\|+\gamma_{21}^{(0)}\|P_{22}^{(0)}\|))I_{n_1},\\
\pi_{22}^{(0)}=P_{22}^{(0)}-\theta(2\gamma_{12}^{(0)}\|P_{12}^{(0)}\|+(\gamma_{12}^{(0)}\|P_{11}^{(0)}\|+\gamma_{21}^{(0)}\|P_{22}^{(0)}\|))I_{n_2},\\
\pi_{12}^{(0)}=P_{12}^{(0)},\quad \pi_{21}^{(0)}=P_{21}^{(0)}
\endgathered
\end{equation*}
and block matrix $P_1=(P_{ij}^{(1)})_{i,j=1,2}$, $(P_{ij}^{(1)})^{\T}=P_{ji}^{(1)}$ with the blocks
\begin{equation}\label{6.1*}
\gathered
P_{11}^{(1)}= e^{-A_{11}^{\T}\theta}(P_{11}^{(0)}
-\theta(P_{12}^{(0)}\widehat A_{21}+\widehat A_{21}^{\T}P_{21}^{(0)}))e^{-A_{11}\theta},
\endgathered
\end{equation}
\begin{equation}\label{6.2*}
\gathered
P_{22}^{(1)}=e^{-A_{22}^{\T}\theta}(P_{22}^{(0)}
-\theta(\widehat A_{12}^{\T}P_{12}^{(0)}+P_{21}^{(0)}\widehat A_{12}))e^{-A_{22}\theta}
\endgathered
\end{equation}
\begin{equation}\label{6.3*}
\gathered
P_{12}^{(1)}= e^{-A_{11}^{\T}\theta}(P_{12}^{(0)}
-\theta(P_{11}^{(0)}\widehat A_{12}+\widehat A_{21}^{\T}P_{22}^{(0)})))e^{-A_{22}\theta}.
\endgathered
\end{equation}


The following proposition is a direct consequence of Theorem 4.1.

{\bf Proposition 6.1.} Let $\Pi_0=(\pi_{ij}^{(0)})_{i,j=1,2}$
be positive definite and 
$Q:=\frac{\Theta_0(\theta)}{\lambda_{\min}(\Pi_0)}+\ln\lambda_{\max}(P_1^{-1}B^{\T}P_0B)<0$,
then \eqref{3.1} is asymptotically stable.

{\it Example.}
\begin{equation}\label{6.4*}
\gathered
\dot x_1(t)=0.01x_1(t)-0.2\sin^2\frac{2\pi t}{\theta} x_2(t),\\
\dot x_2(t)=-0.1x_2(t)+0.2\cos^2\frac{2\pi t}{\theta}x_1(t)
\endgathered
\end{equation}
Choose $\theta=0.09$, $P_{11}^{(0)}=18$, $P_{12}^{(0)}=P_{21}^{(0)}=-7$, $P_{22}^{(0)}=12$,
then $\mu_1=0.01$, $\mu_2=-0.1$, $\delta_1=-0.01$, $\delta_2=0.1$, $\gamma_{12}^{(0)}=\gamma_{21}^{(0)}=0.2$. By the direct calculations we obtain that
\begin{equation*}
\gathered
\Pi_0=\begin{pmatrix}
17.208&-7\\
-7&11.208
\end{pmatrix},\quad
P_1=\begin{pmatrix}
18.0934&-7.0025\\
-7.0024&12.0897
\end{pmatrix}
\endgathered
\end{equation*}
and $\eta_{11}^{(0)}=26.3132$, $\eta_{22}^{(0)}=20.8924$,
$\eta_{12}^{(0)}=7.1037$, $\eta_{21}^{(0)}=7.1036$,
$\alpha_{11}^{(0)}=0.5263$,
$\alpha_{12}^{(0)}=0.0526$,
$\alpha_{21}^{(0)}=0.0418$,
$\alpha_{22}^{(0)}=0.4178$, $Q=-0.00004279<0$.
Therefore, the system is asymptotically stable.
At the same time, the first subsystem is unstable, which makes it impossible to apply a small-gain theorem.

\subsection{Comparison with the small-gain conditions}
Consider \eqref{3.1} with
 $\int_{0}^{\theta}A_{ij}(t)\,dt=0$ for $i\ne j$ and constant dwell-time $T_k=\theta$. Since the Lyapunov vector functions or small-gain results are only applicable when the independent subsystems are asymptotically stable, we assume that
$r_{\sigma}(e^{\theta A_{ii}}B_{ii})<1$ for  $i=1,2$. Now \eqref{6.0}
reduces to two inequalities
\begin{equation}\label{6.5}
\gathered
e^{\theta A_{ii}^{\T}}B_{ii}^{\T}P_{ii}B_{ii}e^{\theta A_{ii}}\prec P_{ii},\quad i=1,2.
\endgathered
\end{equation}
To apply Theorem 4.1, we assume $P_{ii}^{(0)}=P_{ii}$, $i=1,2$
 $P_{12}^{(0)}=0$, then from
\eqref{6.1*}--\eqref{6.3*}
we obtain
\begin{equation*}
\gathered
P_{11}^{(1)}=e^{-\theta A_{11}^{\T}}P_{11}e^{-\theta A_{11}},\quad
P_{22}^{(1)}=e^{-\theta A_{22}^{\T}}P_{22}e^{-\theta A_{22}},\quad P_{12}^{(1)}=0.
\endgathered
\end{equation*}
Let us define
\begin{equation*}
\gathered
\Phi=\begin{pmatrix}
e^{\theta A_{11}}P_{11}^{-1}e^{\theta A_{11}^{\T}}&0\\
0&e^{\theta A_{22}}P_{22}^{-1}e^{\theta A_{22}^{\T}}
\end{pmatrix}
\begin{pmatrix}
B_{11}^{\T}&B_{21}^{\T}\\
B_{12}^{\T}&B_{22}^{\T}
\end{pmatrix}
\begin{pmatrix}
P_{11}&0\\
0&P_{22}
\end{pmatrix}
\begin{pmatrix}
B_{11}&B_{12}\\
B_{21}&B_{22}
\end{pmatrix}
\endgathered
\end{equation*}
A consequence of Theorem 4.1 is the following

{\bf Proposition 6.2.} If $\int_{0}^{\theta}A_{ij}(t)\,dt=0$ for $i\ne j$,  $r_{\sigma}(\Phi)<1$, the conditions of Assumption 4.1 and the following inequalities hold
\begin{equation}\label{StabCond6.1}
\gathered
\theta<\min\Big\{\frac{\lambda_{\min}(P_{11})}{\varrho},\frac{\lambda_{\min}(P_{22})}{\varrho}\Big\},\\
\frac{\Theta_0(\theta)}{\min\{\lambda_{\min}(P_{11})-\theta\varrho,
\lambda_{\min}(P_{22})-\theta\varrho\}}<-\ln\lambda_{\max}(\Phi),
\endgathered
\end{equation}
where $\varrho=\gamma_{12}\|P_{11}\|+\gamma_{21}\|P_{22}\|$.
Then \eqref{3.1} is asymptotically stable.

To compare the obtained results with the results known in the literature, obtained on the basis of the ISS approach or Lyapunov vector function (small-gain conditions), we consider  \eqref{3.1} without impulsive action, i.e. $B_{ii}=I$, $B_{ij}=0$ for $i,j=1,2$, $i\ne j$
and such that
 $\int_{0}^{\theta}A_{ij}(t)\,dt=0$ for $i\ne j$.

Since the  Lyapunov vector function or small-gain results are only applicable when the independent subsystems are asymptotically stable, so we assume, that
$\max\{\Ree\lambda\,|\,\lambda\in\sigma(A_{ii})\}<0$ for $i=1,2$. For given symmetric and positive-definite matrices $Q_i$, $i=1,2$ consider the linear algebraic Lyapunov equations
\begin{equation}\label{6.6}
\gathered
A_{ii}^{\T}P_{ii}+P_{ii}A_{ii}=-Q_i.
\endgathered
\end{equation}
It is known that under our assumptions for $A_{ii}$, $i=1,2$ these equations have unique solutions in the form of symmetric positive-definite matrices $P_{ii}$.

{\it Remark 6.1.} Solutions of matrix algebraic equations \eqref{6.6} satisfy linear matrix inequalities \eqref{6.5}, up to $O(\theta^2)$.

In this case $\Phi=\diag\{e^{\theta A_{11}}P_{11}^{-1}e^{\theta A_{11}^{\T}}P_{11},e^{\theta A_{22}}P_{22}^{-1}e^{\theta A_{22}^{\T}}P_{22}\}$.
Since, we assume that $\theta$ is a sufficiently small positive number, we note that
\begin{equation*}
\gathered
\Phi=I-\theta\diag\{P_{11}^{-1}Q_{1},P_{22}^{-1}Q_{2}\}+O(\theta^2),
\endgathered
\end{equation*}
therefore
\begin{equation*}
\gathered
-\ln\lambda_{\max}(\Phi)=\theta\min\{\lambda_{\min}(P_{11}^{-1}Q_{1}),\lambda_{\min}(P_{22}^{-1}Q_{2})\}+O(\theta^2)
\endgathered
\end{equation*}
 is positive-defined for sufficiently small $\theta>0$. On the other hand, it is easy to show that
$\frac{\Theta_0(\theta)}{\min\{\lambda_{\min}(P_{11})-\theta\varrho,\lambda_{\min}(P_{22})-\theta\varrho\}}=O(\theta^2)$, hence it follows that there exists $\theta^*>0$,
such that for all $\theta\in(0,\theta^*)$ the conditions of
Proposition 6.1 are satisfied. Thus, we come to an important Corollary.

{\bf Corollary 6.1.} System \eqref{3.1} without impulsive action and  $\int_0^{\theta}A_{ij}(t)\,dt=0$ is asymptotically stable for $\theta\in(0,\theta^*)$
if $\theta^*>0$ is small enough.

Note that the $\theta^*$ is determined from the conditions \eqref{StabCond6.1}.

We apply the same Lyapunov functions $V_1(x_1)=x_1^{\T}P_{11}x_1$, $V_2(x_2)=x_2^{\T}P_{22}x_2$ to study the stability of the coupled system \eqref{3.1} without impulsive action, using the small-gain theorem in \cite{EdLiWa2000}. Consider estimates of derivatives
$\dot V_i$ along solutions of \eqref{3.1}. Taking into account Assumption 4.1  and using the Cauchy--Bunyakovsky inequality, we define
\begin{equation}\label{DiffIn}
\gathered
\dot V_i(x_i)=-x_i^{\T}Q_ix_i+2x_i^{\T}P_{ii}A_{ij}(t)x_j\\=
-(P_i^{1/2}x_i)^{\T}P_{ii}^{-1/2}Q_iP_{ii}^{-1/2}P_{ii}^{1/2}x_i+2(P_{ii}^{1/2}x_i)^{\T}P_{ii}^{1/2}A_{ij}(t)P_{jj}^{-1/2}P_{ii}^{1/2}x_j\\
\le-\lambda_{\min}(P_{ii}^{-1/2}Q_iP_{ii}^{-1/2})V_i(x_i)+2\|P_{ii}\|^{1/2}\|P_{jj}\|^{-1/2}\gamma_{ij}V_{ii}^{1/2}(x_i)V_{jj}^{1/2}(x_j)\\
=-\lambda_{\min}(P_{ii}^{-1}Q_i)V_i(x_i)+2\|P_{ii}\|^{1/2}\|P_{jj}\|^{-1/2}\gamma_{ij}V_i^{1/2}(x_i)V_j^{1/2}(x_j).
\endgathered
\end{equation}
Here $i\ne j$, $i,j=1,2$.
To check the conditions of the small-gain theorem from \cite{EdLiWa2000} (Theorem 4), we choose
\begin{equation*}
\gathered
\chi_i(r)=\Big(\frac{2\gamma_{ij}\|P_{ii}\|^{1/2}\|P_{jj}\|^{-1/2}}{\lambda_{\min}(P_{ii}^{-1}Q_i)}+\varepsilon\Big)^2r.
\endgathered
\end{equation*}
Then Theorem 4 from \cite{EdLiWa2000} lead us to the following sufficient conditions for the asymptotic stability of \eqref{3.1} (small-gain condition)
\begin{equation}\label{SmallG}
\gathered
\gamma_{12}\gamma_{21}<\frac{1}{4}\lambda_{\min}(P_{11}^{-1}Q_1)\lambda_{\min}(P_{22}^{-1}Q_2)
\endgathered
\end{equation}
{\it Remark 6.2.} The  method of Lyapunov vector functions lead us to the same stability condition. Indeed,
\eqref{DiffIn} in the new variables
$y_i(t)=V_i^{1/2}(x_i(t))$ leads to a linear system of differential inequalities
\begin{equation*}
\gathered
\dot y_i(t)\le
-\frac{1}{2}\lambda_{\min}(P_{ii}^{-1}Q_i)y_i(t)+\|P_{ii}\|^{1/2}\|P_{jj}\|^{-1/2}\gamma_{ij}y_j(t).
\endgathered
\end{equation*}
 Application of the comparison principle again leads to \eqref{SmallG}. From \eqref{SmallG} follows that the small-gain conditions are not depend on $\theta$. Therefore, it is possible to choose the parameters of the system \eqref{3.1}, such that \eqref{SmallG} is not satisfied, however, based on Corollary 6.1, this system is asymptotically stable for sufficiently small $\theta$.
Thus, we conclude that our approach leads to less conservative stability conditions than the known small-gain conditions.

{\it Example.}
Consider a second-order linear system
\begin{equation}\label{Example1}
\gathered
\dot x_1(t)=-0.2x_1(t)+0.15a_{12}(t) x_2(t),\\
\dot x_2(t)=-0.1x_2(t)+ 0.2a_{21}(t)x_1(t)
\endgathered
\end{equation}
where $a_{ij}\in C(\Bbb R)$, $\|a_{ij}\|_{C[0,\theta]}\leq 1$, $\int\limits_0^{\theta}a_{ij}(t)\,dt=0$, then $\mu_1=-0.2$, $\mu_2=-0.1$, $\delta_1=0.2$, $\delta_2=0.1$ and the small-gain condition
$0.03\|a_{12}\|_{C[0,\theta]}\|a_{21}\|_{C[0,\theta]}<0.02$ is not satisfied.

Choose $\theta=0.5$, $P_{11}^{(0)}=2.5$, $P_{22}^{(0)}=5$,
$P_{12}^{(0)}=0$, $M_1=M_2=N_1=N_2=1$. By a direct calculation $\gamma_{12}^{(0)}=0.15$,
$\gamma_{21}^{(0)}=0.2$,
\begin{equation*}
\gathered
\Pi_0=\begin{pmatrix}
1.8125&0\\
0&4.3125
\end{pmatrix},\quad
P_1=\begin{pmatrix}
3.0535&0\\
0&5.5259
\end{pmatrix}
\endgathered
\end{equation*}
Since $Q=-0.0050178428<0$, the considered system asymptotically stable.


Note that in our example \eqref{Example1}, the functions $a_{ij}$ are actually unknown, we know only restrictions on their norm and mean value. Thus, Theorem 4.1 provides conditions for a whole class of systems \eqref{Example1}. In this sense, Theorem 4.1 is similar to the well-known small-gain theorems for coupled systems.

\section{Numerical examples} 


Here we consider some examples of application of the main results from Sections 4 and 5 to fourth order systems. Example 1 is an impulsive system, with one unstable subsystem. In this case, small-gain conditions are not applicable. Example 2 illustrates the possibility to apply Theorem 4.1 in the case when the continuous and discrete dynamics are both unstable at constant dwell-times, and example 3 is the general case of non-constant dwell-times.

{\it Example 1.} Consider a linear fourth-order impulsive system \eqref{3.1} with the matrices
\begin{equation*}
\gathered
A_{11}=
\begin{pmatrix}
0.1&0.05\\
0.05&0.1
\end{pmatrix},
\quad
A_{22}=
\begin{pmatrix}
-1&0.01\\
0.01&-1
\end{pmatrix}
\endgathered
\end{equation*}
\begin{equation*}
\gathered
A_{12}(t)=-2\sin^2(\omega t)I,\quad
A_{21}(t)=2\sin^2(\omega t)I;
\endgathered
\end{equation*}

\begin{equation*}
\gathered
B_{11}=
\begin{pmatrix}
0.98&0\\
0&0.98
\end{pmatrix},
\quad
B_{22}=
\begin{pmatrix}
1.02&0\\
-0.01&1.01
\end{pmatrix}
\endgathered
\end{equation*}
\begin{equation*}
\gathered
B_{12}=
\begin{pmatrix}
-0.01&0\\
0.02&0
\end{pmatrix},
\quad
B_{21}=
\begin{pmatrix}
0&0.01\\
-0.05&0
\end{pmatrix}.
\endgathered
\end{equation*}
Here, $\tau_{k+1}-\tau_k=\theta=0.2$, $\omega=\frac{2\pi}{\theta}$. To check the asymptotic stability conditions of Theorem 6.1, we choose
$N=50$, $P_{11}^{(0)}=18I$, $P_{12}^{(0)}=-7 I$, $P_{22}^{(0)}=12I$.
Then $\min_{k}\lambda_{\min}(\Pi_k)=7.0322$, $Q=-0.0464<0$. Hence \eqref{3.1} is asymptotically stable.
Note, that the independent subsystem
\begin{equation}\label{subsystems}
\gathered
\dot x_2(t)=A_{22}x_2(t),\quad t\ne\tau_k,\\
x_2(t^+)=B_{22}x_2(t),\quad t=\tau_k,\\
\endgathered
\end{equation}
is not stable due to the fact that $r_{\sigma}(e^{A_{22}\theta}B_{22})>1$.

{\it Example 2.} Consider a linear impulsive system \eqref{3.1} with the matrices
\begin{equation*}
\gathered
A_{11}=
\begin{pmatrix}
-1&0\\
0&-1
\end{pmatrix},
\quad
A_{22}=
\begin{pmatrix}
0.1&0\\
0&0.1
\end{pmatrix}
\endgathered
\end{equation*}
\begin{equation*}
\gathered
A_{12}(t)=
\begin{pmatrix}
0.2\cos(\omega t)&-0.2\sin(\omega t)\\
0.2\sin(\omega t)&0.2\cos(\omega t)
\end{pmatrix},
\\
A_{21}(t)=
\begin{pmatrix}
0.1\cos(\omega t)&-0.1\sin(\omega t)\\
0.1\sin(\omega t)&0.1\cos(\omega t)
\end{pmatrix}
\endgathered
\end{equation*}
\begin{equation*}
\gathered
B_{11}=
\begin{pmatrix}
1.2&0.1\\
-0.1&1.5
\end{pmatrix},
\quad
B_{22}=
\begin{pmatrix}
0.5&0.05\\
-0.05&-0.5
\end{pmatrix}
\endgathered
\end{equation*}
\begin{equation*}
\gathered
B_{12}=
\begin{pmatrix}
0.04&0.1\\
0.1&0.04
\end{pmatrix},
\quad
B_{21}=
\begin{pmatrix}
0.05&0.1\\
0.2&0.05
\end{pmatrix}
\endgathered
\end{equation*}
Here, $\tau_{k+1}-\tau_k=\theta=0.5$, $\omega=\frac{2\pi}{\theta}$. To check the asymptotic stability conditions, obtained in Theorem 6.1, we choose
$N=3$, $P_{11}^{(0)}=I$, $P_{12}^{(0)}=0$, $P_{22}^{(0)}=I$.
Then,
\begin{equation*}
\gathered
P_{3}=
\begin{pmatrix}
2.7172&0&-0.0008&-0.0207\\
0&2.7172&0.0207&-0.0008\\
-0.0008&0.0207&0.902&0\\
-0.0207&-0.0008&0&0.902
\end{pmatrix},
\endgathered
\end{equation*}
$\min_{k}\lambda_{\min}(\Pi_k)=0.7793$, $Q=-0.0383<0$. Therefore, the linear impulsive system \eqref{3.1} is asymptotically stable.
Consider separately the continuous dynamics of a linear impulsive system, which is described by a linear time-variant ODEs
\begin{equation}\label{example2}
\gathered
\dot x_1(t)=A_{11}x_1(t)+A_{12}(t)x_2(t),\\
\dot x_2(t)=A_{21}(t)x_1(t)+A_{22}x_2(t),
\endgathered
\end{equation}
We denote
\begin{equation*}
U(\omega t)=
\begin{pmatrix}
\cos(\omega t)&-\sin(\omega t)\\
\sin(\omega t)&\cos(\omega t)
\end{pmatrix}
\end{equation*}
and rewrite \eqref{example2} as
\begin{equation}\label{example21}
\gathered
\dot x_1(t)=-x_1(t)+0.2U(\omega t)x_2(t),\\
\dot x_2(t)=0.1U(\omega t)x_1(t)+0.1x_2(t),
\endgathered
\end{equation}
Consider the Lyapunov--Chetaev function $v(x_1,x_2)=2\|x_2\|^2-\|x_1\|^2$, the total derivative of which is
\begin{equation*}
\gathered
\dot v(x_1,x_2)=2(0.2 x_2\|^2+\|x_1\|^2+0.4\sin(\omega t)x_1^{\T}Jx_2)\\ \ge
2(0.2\|x_2\|^2+\|x_1\|^2-0.4\|x_1\|\|x_2\|)>0,\quad \text{for all }(x_1,x_2)\ne (0,0).
\endgathered
\end{equation*}
Here, $
J=\begin{pmatrix}
0&1\\
-1&0
\end{pmatrix}
$ is the symplectic unit.
Hence, the ODE \eqref{example21} is unstable.
The discrete dynamics is unstable since $r_{\sigma}(B)=1.4722$.

{\it Example 3.} Consider an example, which illustrate the application of Theorem 5.1 to the case of a non-periodic impulsive system. Let
\begin{equation*}
\gathered
A_{11}=
\begin{pmatrix}
-1&0\\
0&-1
\end{pmatrix},
\quad
A_{22}=
\begin{pmatrix}
0.1&0\\
0&0.1
\end{pmatrix}
\endgathered
\end{equation*}
\begin{equation*}
\gathered
A_{12}(t)=
\begin{pmatrix}
0.05\cos(\omega t)&-0.05\sin(\omega t)\\
0.05\sin(\omega t)&0.05\cos(\omega t)
\end{pmatrix},
\\
A_{21}(t)=
\begin{pmatrix}
0.05\cos(\omega t)&-0.05\sin(\omega t)\\
0.05\sin(\omega t)&0.05\cos(\omega t)
\end{pmatrix}
\endgathered
\end{equation*}
\begin{equation*}
\gathered
B_{11}=
\begin{pmatrix}
1.2&0.1\\
-0.1&1.5
\end{pmatrix},
\quad
B_{22}=
\begin{pmatrix}
0.25&0.05\\
-0.05&-0.25
\end{pmatrix}
\endgathered
\end{equation*}
\begin{equation*}
\gathered
B_{12}=
\begin{pmatrix}
0.04&0.1\\
0.1&0.04
\end{pmatrix},
\quad
B_{21}=
\begin{pmatrix}
0.05&0.1\\
0.2&0.05
\end{pmatrix}
\endgathered
\end{equation*}
For $\theta=1$, $\theta_1=0.8$, $\theta_2=1.2$, we choose
$N=7$, $P_{11}^{(0)}=I$, $P_{12}^{(0)}=0$, $P_{22}^{(0)}=I$. Then, $N_3=4$, $N_4=8$, $\min_{m=\overline{0,11}}\lambda_{\min}(\Pi_m^{(l)})=0.7585>0$,
$\max_{l=\overline{0,9}, m=\overline{6,12}}Q_{m}^{(l)}=-0.1399<0$. Thus, by Theorem 5.1, the linear system \eqref{3.1} is asymptotically stable only if the dwell-times $T_k=\tau_{k+1}-\tau_k$ belong to the interval $[\theta_1,\theta_2]$. We note that in this example both continuous and discrete dynamics are not stable. The instability of continuous dynamics is proved using the Lyapunov--Chetaev function $v(x_1,x_2)=2\|x_2\|^2-\|x_1\|^2$ as done above. The instability of discrete dynamics follows from the fact that $r_{\sigma}(B)=1.4746$.

\section{Discussion}
Theorems 4.1 and 5.1 are the main result of the paper and establish sufficient conditions of the exponential stability of the linear time-variant impulsive system \eqref{3.1}. The proposed approach to the study of a coupled impulsive system significantly expands the capabilities of the method of Lyapunov vector functions and allows one to study the asymptotic stability of a linear impulsive system with unstable subsystems. Given examples in Section 6 show that our results significantly expand the known methods for studying impulsive systems developed in \cite{dash-miron,hadad,liu,dv-sl-1}. In addition, the obtained results are applied in the case when the continuous and discrete dynamics of an impulsive system are both unstable. In the case when dwell-times are non constant for studying the stability of a linear impulsive system, the classical Floquet theory turns out to be inapplicable.
These results significantly expand the possibilities of the direct Lyapunov method in the context of theorems from \cite{sam-per}.

In this case, Theorem 5.1 allows us to deduce the Lyapunov asymptotic stability.
It is of interest to generalize the proposed approaches for the construction of Lyapunov functions for nonlinear coupled systems, as well as when expanding the assumptions about the number of independent subsystems.

\section{Appendix}

{\it Proof} of Lemma 4.1.
Using the expressions for $P_{ij}^{(m)}$ we obtain 
\begin{equation*}
\gathered
v(t,x_1,x_2)
=z_{1m}^{\T}\Big(P_{11}^{(m)}-\int\limits_{mh}^t(P_{12}^{(m)}A_{21}(s)+A_{21}^{\T}(s)P_{21}^{(m)})\,ds\Big)z_{1m}\\
+2z_{1m}^{\T}\Big(P_{12}^{(m)}-\int\limits_{mh}^t(P_{11}^{(m)}A_{12}(s)+A_{21}^{\T}(s)P_{22}^{(m)})\,ds\Big)z_{2m}\\
+z_{2m}^{\T}\Big(P_{22}^{(m)}-\int\limits_{mh}^t(A_{12}^{\T}(s)P_{12}^{(m)}+P_{21}^{(m)}A_{12}(s))\,ds\Big)z_{2m}.
\endgathered
\end{equation*}
Applying the Cauchy-Bunyakovsky inequality, we obtain
\begin{equation*}
\gathered
\Big|z_{1m}^{\T}\int\limits_{mh}^t(P_{12}^{(m)}A_{21}(s)+A_{21}^{\T}(s)P_{21}^{(m)})\,dsz_{1m}\Big|
\le 2\|P_{12}^{(m)}\|\gamma_{21}^{(m)}\|z_{1m}\|^2,
\endgathered
\end{equation*}
\begin{equation*}
\gathered
\Big|z_{1m}^{\T}\int\limits_{mh}^t(P_{11}^{(m)}A_{12}(s)+A_{21}^{\T}(s)P_{22}^{(m)})\,dsz_{2m}\Big|
\le (\|P_{11}^{(m)}\|\gamma_{12}^{(m)}+\|P_{22}^{(m)}\|\gamma_{21}^{(m)})\|z_{1m}\|\|z_{2m}\|,
\endgathered
\end{equation*}
\begin{equation*}
\gathered
\Big|z_{2m}^{\T}\int\limits_{mh}^t(A_{12}^{\T}(s)P_{12}^{(m)}+P_{21}^{(m)}A_{12}(s))\,dsz_{2m}\Big|
\le 2\|P_{12}^{(m)}\|\gamma_{12}^{(m)}\|z_{2m}\|^2
\endgathered
\end{equation*}
Thus,
\begin{equation*}
\gathered
v(t,x_1,x_2)\ge z_{1m}^{\T}P_{11}^{(m)}z_{1m}+2z_{1m}^{\T}P_{12}^{(m)}z_{2m}+z_{2m}^{\T}P_{22}^{(m)}z_{2m}
-2\|P_{12}^{(m)}\|\gamma_{21}^{(m)}\|z_{1m}\|^2\\
-2(\|P_{11}^{(m)}\|\gamma_{12}^{(m)}+\|P_{22}^{(m)}\|\gamma_{21}^{(m)})\|z_{1m}\|\|z_{2m}\|-
2\|P_{12}^{(m)}\|\gamma_{12}^{(m)}\|z_{2m}\|^2.
\endgathered
\end{equation*}
Applying the Cauchy inequality, we get
\begin{equation*}
\gathered
v(t,x_1,x_2)\ge z_{1m}^{\T}P_{11}^{(m)}z_{1m}+2z_{1m}^{\T}P_{12}^{(m)}z_{2m}+z_{2m}^{\T}P_{22}^{(m)}z_{2m}\\
-(2\|P_{12}^{(m)}\|\gamma_{21}^{(m)}+(\|P_{11}^{(m)}\|\gamma_{12}^{(m)}+\|P_{22}^{(m)}\|\gamma_{21}^{(m)}))\|z_{1m}\|^2\\
-(2\|P_{12}^{(m)}\|\gamma_{12}^{(m)}+(\|P_{11}^{(m)}\|\gamma_{12}^{(m)}+\|P_{22}^{(m)}\|\gamma_{21}^{(m)}))\|z_{2m}\|^2
\ge\lm_{\min}(\Pi_m)\|z_m\|^2.
\endgathered
\end{equation*}
The upper bound for the Lyapunov function $v(t,x_1,x_2)$ is proved similarly. 

\bibliography{mybibfile}

\end{document}